\documentclass [a4paper,12pt]{amsart}
\usepackage{graphicx}
\usepackage[ansinew]{inputenc}    
\usepackage[all]{xy}
\usepackage{multicol}
\usepackage{stackengine}
\usepackage{tikz-cd}

\usepackage{mathabx}
\usepackage{amssymb,amscd}
\usepackage{mathrsfs}

\usepackage{enumitem}
\usepackage{verbatim}
\usepackage{hyperref}

\newtheorem{thm}{Theorem}[section]
\newtheorem{defn}[thm]{Definition}

\newtheorem{exam}[thm]{Example}
\def\im{\operatorname{Im}}

\def\Fix{\operatorname{Fix}}

\def\C{\mathbb C}
\def\R{\mathbb R}

\def\Z{\mathbb Z}

\def\dim{\operatorname{dim}}

\def\id{\operatorname{id}}

\usepackage{graphicx}
\usepackage[ansinew]{inputenc}    
\usepackage[all]{xy}
\usepackage{amssymb,amscd}

\newtheorem{cor}[thm]{Corollary}

\newtheorem{lem}[thm]{Lemma}
\newtheorem{prop}[thm]{Proposition}

\theoremstyle{definition}
\newtheorem{ex}[thm]{Example}
\newtheorem{defi}[thm]{Definition}

\newtheorem{rem}[thm]{Remark}

\def\C{\mathbb C}
\def\R{\mathbb R}

\def\Z{\mathbb Z}

\def\dim{\operatorname{dim}}

\def\id{\operatorname{id}}

\usepackage{color}

\thanks{The first author has been partially supported by CAPES. The second author has been suported by Grant PGC2018-094889-B-100 funded by MCIN/AEI/ 10.13039/501100011033 and by ``ERDF A way of making Europe''. The third has been supported by FAPESP-Grant 2022/15458-1. The third and fourth author have been partially supported by  FAPESP Grant 2019/07316-0.}
\keywords{Reflection groups, singular mappings, multiple points}
\subjclass[2000]{Primary 32S25; Secondary 58K40, 32S50}

 \everymath{\displaystyle}
\begin{document}
\title{Double points and image of reflection maps}

\author{J. R. Borges-Zampiva, B. Or\'efice-Okamoto, G. Pe\~nafort Sanchis,
	 J.N. Tomazella}

\address{Departamento de Matem\' atica, Universidade Federal de S\~ao Carlos, Caixa Postal 676, 13560-905, S\~ao Carlos, SP, BRAZIL}
\email{jrbzampiva@estudante.ufscar.br}

\address{Departament de Matem\`atiques,
	Universitat de Val\`encia, Campus de Burjassot, 46100 Burjassot
	SPAIN.}
\email{guillermo.penafort@uv.es}

\address{Departamento de Matem\'atica, Universidade Federal de S\~ao Carlos, Caixa Postal 676,
	13560-905, S\~ao Carlos, SP, BRAZIL}
\email{brunaorefice@ufscar.br}

\address{Departamento de Matem\'atica, Universidade Federal de S\~ao Carlos, Caixa Postal 676,
	13560-905, S\~ao Carlos, SP, BRAZIL}
\email{jntomazella@ufscar.br}

\maketitle

\begin{abstract}
A reflection mapping is a singular holomorphic mapping obtained by restricting the quotient mapping of a complex reflection group. We study the analytic structure of double point spaces of reflection mappings. In the case where the image is a hypersurface, we obtain explicit equations for the double point space and for the image as well. In the case of surfaces in $\C^3$, this gives a very efficient method to compute the Milnor number and delta invariant of the double point curve.

\end{abstract}

\section{Introduction}
In the theory of singular mappings, there are few known examples which are degenerate and also have desirable properties. The problem, rather than an actual lack of mappings of this kind, seems to be the difficulty of the calculations.
Reflection mappings were introduced in the prequel of this work  \cite{Penafort-Sanchis:2016a}, as a means to produce degenerate mappings which are easy to understand. This stablished the existence of the examples we were looking for in certain dimensions. In this paper we go one step further and show how to compute the double point spaces---and, in the hypersurface case, the image as well--- for reflection mappings.  In short, a reflection mapping is a holomorphic singular mapping $\mathcal Y\to \C^p$ obtained by restricting the quotient map of a reflection group $W$ to a submanifold $\mathcal Y$ of the vector space $\mathcal V$ where the group acts (see Section \ref{secPreliminaries} for details).  The idea is that the group action, together with the way $\mathcal Y$ sits in $\mathcal V$, must encode the geometry of the mapping. 

The reflection mapping class contains the first families of quasi-homogeneous finitely determined map-germs $(\C^n,0)\to(\C^{p},0)$ having unbounded multiplicity, for arbitrary $n$ and $p=2n-1$ or $p=2n$ \cite[Theorems 9.5 and 9.6]{Penafort-Sanchis:2016a} (see Example \ref{exRefMapsUnbounded} for low dimensional examples). These mappings where later shown by Ruas and Silva to be counterexamples to a long standing conjecture of Ruas \cite{Ruas:2017}. Brasselet, Thuy and Ruas used them to show the density of the finitely determined mappings among certain spaces of quasi-homogeneous mappings \cite{BrNgRu2019}. Silva noticed that these reflection mappings show that the topological type of a generic transverse slice, in the case of $\mathcal A$-finite quasihomogeneous mappings, is not determined by their weights and degrees \cite{MR4706774}. Rodrigues Hernandes and Ruas have shown that every finitely determined monomial  mapping $(\C^n,0)\to(\C^p,0)$, with $p\geq 2n$, is a reflection mapping \cite{MR3975520}. In contrast, for $p<2n-1$, the only finitely determined reflection mappings are folding maps \cite[Theorem 8.5]{Penafort-Sanchis:2016a}.  Reflection mappings were also shown to satisfy an extended version of L\^e's conjecture, a rather misterious problem relating injectivity and corank \cite[Proposition 4.3]{Penafort-Sanchis:2016a}. 

Apart from the theory of singular mappings, reflection mappings have proven to be relevant in differential geometry. The original idea is due to Bruce and Wilkinson \cite{MR757467, MR1129024}, who noticed that the singularities produced by folding a surface with respect to a plane (in our terminology these are $\Z/2$-reflection mappings) reveal interesting extrinsic information of the surface. This idea has since been extended to the other spaces \cite{Barajas2020,MR2768804,MR3847079} as well as to  different reflection groups \cite{PeTa2020}. The reflection mappings used by Bruce and Wilkinson are known as folding mappings and have been widely studied (apart from the previous references, see \cite{MR1646521}). They belong to a subclass called reflected graphs \cite[Definition 13]{Penafort-Sanchis:2016a}, which we describe in Example \ref{exReflectedGraphs}. This class includes other previously studied classes of mappings, such as double folds \cite{MR2431646,MR3300299} and $k$-folding mappings \cite{PeTa2020}. Double folds provided the first examples of finitely determined quasi-homogeneous corank two map-germs. Inspecting the classification  of all simple map germs $(\C^2,0)\to(\C^3,0)$ \cite{MR772717} reveals that they are either folding mappings or $3$-folding mappings.

Now we can discuss the content of the present work. Singular mappings are understood via the study of their multiple point spaces, the most fundamental of which are the double point spaces. Moreover, given the coordinate functions of a map-germ $(\C^n,0)\to(\C^{n+1},0)$, one would like to know the equation of its image. As it turns out, both the image and the double points spaces become degenerate very fast as the complexity of the coordinate functions increases. We show how to study these objects for reflection mappings, in a much simpler way than in the case of arbitrary singular mappings. In the hypersurface case, the image of reflection mappings is described as well.

 For reflection mappings, there are three decompositions of double point spaces indexed by  the reflection group, namely
 \[K_2(f)=\bigcup_{\sigma\in W\setminus\{1\}}K_2^\sigma(f),\qquad D^2(f)=\bigcup_{\sigma\in W\setminus\{1\}}D_2^\sigma(f),
\qquad D(f)=\bigcup_{\sigma\in W\setminus\{1\}}D_\sigma(f).\]
These decompositions appeared in \cite[Sections 6 and 7]{Penafort-Sanchis:2016a} already, but just at the set theoretical level. In contrast, here the analytic structure of the branches $K_2^\sigma(f)$, $D_2^\sigma(f)$ and $D_\sigma(f)$ is described explicitly. These analytic structures are a fundamental part of the theory of singular mappings, which requires $K_2(f)$, $D^2(f)$ and $D(f)$ to be complex spaces, sometimes with a non-reduced structure indicating higher degeneracy. 

For arbitrary mappings $\mathcal Y\to \C^{n+1}$ with $\dim \mathcal Y=n$, the double point space $D(f)\subset \mathcal Y$ and the image $\im f\subset \C^{n+1}$ are hypersurfaces (or more precisely, complex spaces defined locally by principal ideals, as they are not always reduced), but computing their equations is usually very hard. One of the main achievements of this work are the explicit formulas
\[\im f=V\Big(\prod_{\sigma\in W}(\sigma L)\circ s\Big),\qquad D(f)=V\Big(\prod_{\sigma\in W\setminus\{1\}}\lambda_\sigma\Big)\]
for reflection mappings. As a bonus, the formula for $D(f)$ gives a much faster way to compute the Milnor number and delta invariant of $D(f)$ for germs of reflection mappings $\mathcal Y^2\to\C^3$.

\section{Notation and preliminaries}\label{secPreliminaries}
In this section we summarize the prerequisites and set the notation for reflection groups and reflection mappings. For a more detailed account and proper citations, we refer to \cite{Penafort-Sanchis:2016a}.
To avoid repetition, we fix the meaning of the symbols and notations summarized in this section, which will not be reintroduced.

\subsection*{Reflection groups and the orbit mapping}Throughout the text, $W$ stands for a (complex) reflection group acting on a $\C$-vector space $\mathcal V$. We  write $p=\dim \mathcal V$.

 The orbit of a subset $S\subseteq \mathcal V$ (or better, the union of the orbits of the points in $S$) is denoted by $WS\subseteq \mathcal V$. 
We adopt the convention that the action of $\sigma\in W$ on a function $H$ in  $\mathcal O_{\mathcal V}$ is
$(\sigma H)(u)=H(\sigma^{-1}u).$
Hence, a zeroset $S=V(J)$ is transformed by $\sigma$ into the set $\sigma S=\{\sigma u\mid u\in S\}=V(\sigma J)$.
The \emph{stabilizer} of $S$ is
\[W^S=\{\sigma\in W\mid \sigma S=S\}.\]
and the \emph{pointwise stabilizer} of $S$  is
\[W_S=\{\sigma\in W\mid \sigma u=u,\text{ for all }u\in S\}.\]
$W^S$ is the maximal subgroup of $W$ acting on $S$. The pointwise stabilizer $W_S$ is a normal subgroup of $W^S$ and the quotient $W^S/W_S$ acts faithfully on $S$.

The union of the reflecting hyperplanes of all reflections in $W$ is called the \emph{hyperplane arrangement} and is written as
$\mathscr A\subset\mathcal V$.
The reflecting hyperplanes induce a partition of $\mathcal V$ into subsets, called facets, consisting of those points contained in exactly the same hyperplanes. The set of facets is denoted by $\mathscr C$ and called the \emph{complex} of $W$.
Clearly, facets $C\in \mathscr C$ are open subsets of linear subspaces, hence have the same tangent everywhere, which may be identified with the closure $\overline C$. Similarly, we often write $C^\bot$ instead of $T_yC^\bot$, for any $y\in C$.

The celebrated Shephard-Todd-Chevalley Theorem characterizes reflection groups as the only subgroups of the group of unitary transformations of $\mathcal V$ for which the quotient mapping $\mathcal V\to \mathcal V/W$ can be realized as a polynomial mapping 
\[\omega\colon \mathcal V\to \C^p\]
whose coordinate functions $\omega_1,\dots,\omega_p$ are homogeneous polynomials (they are a set of generators of the ring of $W$-invariant polynomials). More geometrically, the map $\omega$ identifies an orbit to a point, that is, for any $u\in\mathcal V$, we have that $\omega^{-1}(\omega(u))=Wu$.

Reflection groups act faithfully, hence $\omega$ is generically $|W|$-to-one, and it is well known that the ramification locus is precisely the hyperplane arrangement. More precissely, given a point $u\in \mathcal V$, contained in a facet $C\in \mathscr C$, we have that $\ker \operatorname{d}\!\omega_u=C^\bot$.

 By the universal property of quotient mappings, the mapping $\omega$ is well defined up to $\mathcal L$-equivalence (that is, up to changes of coordinates in the target). As it will become clear, replacing $\omega$ by an $\mathcal L$-equivalent map does not change the $\mathcal A$-class of the associated reflection mappings. Since we study reflection mappings up to $\mathcal A$-equivalence, we may pretend $\omega$ to be unique, and abusively call it \emph{the orbit map} of $W$.
 
  	\begin{ex}[Products of cyclic groups]\label{exProductOfCyclicGroups}
 		The product $\Z/{d_1}\times\dots\times \Z/{d_p}$ is a reflection group acting on $\C^p$ by \[(a_1,\dots, a_p)\cdot (u_1,\dots,u_p)=(\xi_1^{a_1}u_1,\dots,\xi_p^{a_p}u_p)\text{, with }\xi_j= e^{\frac{2\pi i}{d_j}}.\] An element $(a_1,\dots, a_p)$ is a reflection if and only if exactly one $a_i$ is nonzero. The orbit map is defined by $\omega(u_1,\dots, u_p)=(u_1^{d_1},\dots, u_p^{d_p}).$
 	\end{ex}

 
 \begin{ex}[Dihedral groups, $D_{2n}$]\label{exD2n}
These are the groups of symmetries of regular polygon of $n$ sides. We will use the group $D_8$, consisting of the identity,  four reflections
%
 \[
 \sigma_1=\left(\begin{array}{cc}1 & 0 \\
 	0 & -1\end{array}\right),\quad
 \sigma_2=\left(\begin{array}{cc}0 & 1 \\ 
 	1 & 0\end{array}\right),\quad
 \sigma_3=\left(\begin{array}{cc}-1 & 0 \\
 	0 & 1\end{array}\right),\quad
\sigma_4=\left(\begin{array}{cc}0 & -1 \\
 	-1 & 0\end{array}\right)
 \]
 and three rotations
 \[
 \rho_1=\left(\begin{array}{cc}0 & -1 \\ 
 	1 & 0\end{array}\right),\quad
 \rho_2=\left(\begin{array}{cc}-1 & 0 \\
 	0 & -1\end{array}\right),\quad
 \rho_3=\left(\begin{array}{cc}0 & 1 \\ 
 	-1 & 0\end{array}\right).
 \]	
 The orbit map of $D_8$ is the mapping $\omega\colon \C^2\to\C^2$, given by
 \[(u,v)\mapsto(u^2+v^2,u^2v^2).\]

 \end{ex}
 
 	\begin{exam}[The group $\mathfrak S_4$ of symmetries of the tetrahedron]\label{examS4Group}
 		Consider a regular tetrahedron centered at the origin of $\R^3$, for example the one with vertices
\[V_1=\left(1, 0, \frac{-1}{\sqrt{2}}\right)\quad V_2=\left(-1, 0, \frac{-1}{\sqrt{2}}\right),\quad V_3=\left( 0,1, \frac{1}{\sqrt{2}}\right),\quad V_4=\left(0,-1, \frac{1}{\sqrt{2}}\right).\]
The group of unitary automorphisms of the tetrahedron is a reflection group, isomorphic as an abstract group to the permutation group $\mathfrak S_4$ of its four vertices. This group is generated by the permutations $(i\;i+1)$, which in matrix form are
 		\[(1\;2)=\left(\begin{array}{ccc}-1 & 0 & 0 \\0 & 1 & 0 \\0 & 0 & 1\end{array}\right),\;
 		(2\;3)=\frac{1}{2}\left(\begin{array}{ccc}1 & -1 & -\sqrt{2} \\-1 & 1 & -\sqrt{2} \\-\sqrt{2} & -\sqrt{2} & 0\end{array}\right),\;
 		(3\;4)=\left(\begin{array}{ccc}1 & 0 & 0 \\0 & -1 & 0 \\0 & 0 & 1\end{array}\right).
 		\]
The group $\mathfrak S_4$ contains $6$ reflections, which are precisely the permutations $(i\;j)$. The orbit map is \[\omega(u,v,w)=\big(u^2+v^2+w^2,(u+v)(u-v)w,(2u^2-w^2)(2v^2-w^2)\big).\]


 	\end{exam}
%
%

\subsection*{Reflection mappings} Take an embedding $h\colon \mathcal X\hookrightarrow \mathcal V$ of an $n$-dimensional complex manifold $\mathcal X$ into $\mathcal V$.  The image of $h$ is an $n$-dimensional complex submanifold 
\[\mathcal Y=h(\mathcal X)\subseteq \mathcal V.\] A \emph{reflection mapping} is a map obtained as the composition of the orbit map $\omega$ and the embedding $h$, that is, 
\[f=\omega\circ h\colon \mathcal X\to \C^p.\]
It is often convenient to replace $h$ by the inclusion $\mathcal Y\hookrightarrow \mathcal V$, obtaining a reflection mapping which, abusively, is also denoted
\[f=\omega\vert_{\mathcal Y}\colon \mathcal Y\to \C^p.\]
The choice between these two equivalent settings will be clear from the context.
To finish fixing our notation, locally at any point, $\mathcal Y$ is defined in $\mathcal V$ by a collection of regular equations, which we write as 
\[L=(L_1,\dots,L_{p-n})=0.\]
(technically, it would be better to express our results in terms of the ideal sheaf $I(\mathcal Y)$ of holomorphic functions vanishing on $\mathcal Y$. For the sake of clarity, we have chosen to ignore this issue, which can be fixed by a standard glueing process).
\begin{ex}[Reflected graphs]\label{exReflectedGraphs}
Take a reflection group $W$ acting on $\mathcal V$, and a mapping $H\colon \mathcal V\to \C^p$. We may regard $W$ as a reflection group acting on $\mathcal V\times \C^p$, by extending the action trivially on $\C^p$, and take $\mathcal Y\subseteq \mathcal V\times \C^p$ to be the graph of $H$ (or equivalently, take the graph embedding $h\colon \mathcal V\to \mathcal V\times \C^p$)  . The resulting reflection mapping $\mathcal V\to \C^n\times \C^p$ is called a \emph{reflected graph} and has the form
\[x\mapsto \big(\omega(x),H(x)\big).\]

A typical and much studied  example of reflected graphs are folding maps \[(x,y)\mapsto (x,y^2, H(x,y)),\quad x\in \C^{n-1},y\in \C.\] Observe that, after a target change of coordinates, $H(x,y)$ can be taken to be of the form $yP(x,y^2)$. Similarly, double folds are $\mathbb Z/2\times \mathbb Z/2$-reflected graphs of the form $(x,y)\mapsto (x^2,y^2,xP_1(x^2,y^2)+yP_2(x^2,y^2)+xyP_3(x^2,y^2))$ and $k$-folding  mappings are $\Z/k$-reflected graphs $(x,y)\mapsto(x,y^k,H(x,y))$ (See the Introduction for references for folding maps, double folds and $k$-folding mappings).

%
%
An interesting example is the $\mathbb Z/k\times \mathbb Z/k$-reflection graph $(x,y)\mapsto (x^k,y^k, xy)$, defined by  the graph of the function $H(x,y)=xy$. It parametrizes the $A_{k-1}$-singularity $\{Z^k=XY\}$, but it does so in a generically $k$-to-one way, as explained in Example \ref{exDegree} (observe that  $A_{k-1}$ being normal prevents it from being parametrized in a generically one-to-one way).

Also, taking the group $D_8$ (Example \ref{exD2n}) and the functions $H_1(x,y)=x+2y$ and $H_2(x,y)=x^2+xy-y^2+x^3+x^2y-2xy^2-y^3$, we obtain the reflection graphs
\[f^{D_8}_1\colon(x,y)\mapsto (x^2+y^2,x^2y^2,x+2y),\]
\[f^{ D_8 }_2\colon (x,y)\mapsto (x^2+y^2,x^2y^2,2x^2+3xy-y^2+2x^3+8x^2y-2xy^2-2y^3).\]
The last two examples are depicted in Figure \ref{imaged82}.

\begin{center}
 \begin{figure}[h]
\includegraphics[scale=0.9]{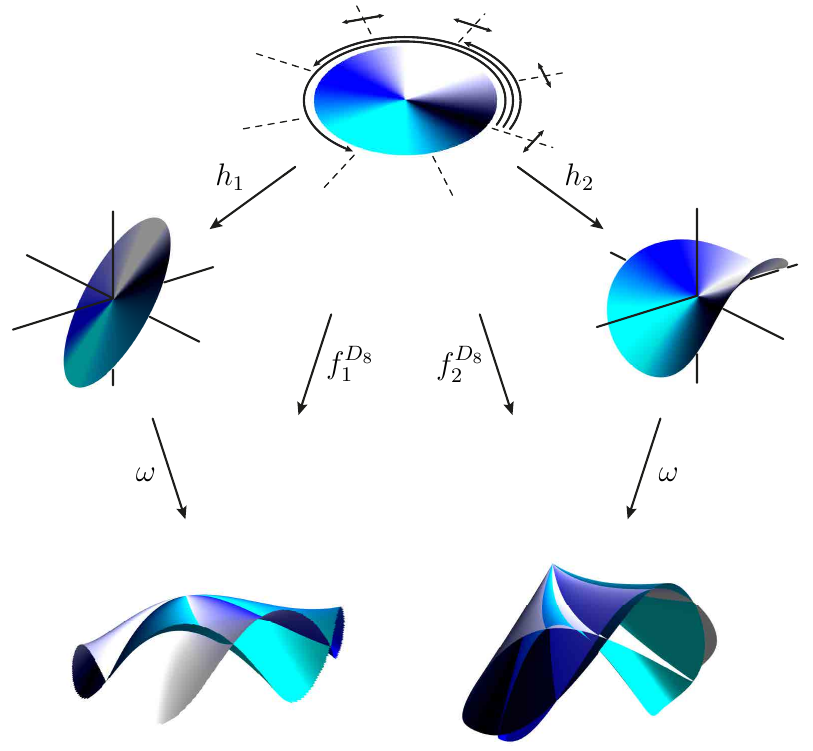}
	\caption{Construction of the reflected graphs $f_1^{D_8}$ and $f_2^{D_8}$.}
	\label{imaged82}
	\end{figure}
\end{center}

\end{ex}

\begin{ex}[Reflection mappings with unbounded multiplicity]\label{exRefMapsUnbounded}
The map-germs
\[f^{(d_1,d_2,d_3)}\colon(x,y)\mapsto (x^{d_1}, y^{d_2},(x+y)^{d_3}),\]
\[f^{(d_1,\dots,d_5)}\colon(x,y,z)\mapsto (x^{d_1}, y^{d_2},z^{d_3},(x+y+z)^{d_4},(x-y+2z)^{d_5}),\]
with $d_i$ pairwise coprime positive integers, belong to a family of map-germs $(\C^n,0)\to(\C^{2n-1},0)$, introduced in  \cite[Theorems 9.5 and 9.6]{Penafort-Sanchis:2016a}), which, to this day, are the only known family of $\mathcal A$-finite map-germs whose coordinate functions have unbounded order. These are the germs studied by Brasselet, Ruas, Silva and Thuy mentioned in the Introduction.

\end{ex}

\subsection*{Unfoldings and $W$-unfoldings}

In the theory of singularities of mappings, the notion of deformation of a space-germ  is replaced by that of unfolding of  a map-germ. An \emph{unfolding} of a map-germ $f\colon (\mathcal Y,S)\to  (\mathcal Z,z)$ is a map-germ \[F\colon (\mathcal Y,S)\times (\C^r,0)\to  (\mathcal Z,z)\times (\C^r,0)\] 
of the form $F(x,t)=(f_t(x),t)$, satisfying $f_0=f$.  Fixed a small representative of $F$ and a value $\epsilon$ of the parameter $t$, the map $f_\epsilon$ is called a \emph{perturbation} of $f$.

A reflection mapping may be perturbed into a non reflection mapping, for example, by perturbing $\omega$ in a way unrelated to the action of the group $W$. At the same time, perturbing $\mathcal Y$ while keeping $\omega$ intact gives a family of reflection mappings containing the original reflection mapping defined by $\mathcal Y$. If we want to study unfoldings of reflection mappings without leaving the reflection mapping seting, these are the deformations we want to consider. In the following lines we formalize this construction.

Let $f\colon \mathcal Y\to \C^p$ be a reflection mapping, and consider a complex submanifold $\widetilde{\mathcal Y}\subseteq \mathcal V\times \Delta$, such that $\widetilde{\mathcal Y}\to \Delta$ defines a trivial deformation of $\mathcal Y$ over an open subset $\Delta\subseteq \C^r$ containing the origin. We may extend the action of $W$ to $\mathcal V\times \C^r$, trivially on the $\C^r$ coordinates, so that the corresponding orbit mapping is $\widetilde\omega=\omega\times \id_{\C^r}$. 
 \begin{defn}
In the above setting, the reflection mapping
\[F=\widetilde\omega_{\vert_{\widetilde{\mathcal Y}}}\colon \widetilde{\mathcal Y}\to \C^p\times \Delta\]
is called a \emph{$W$-unfolding} of $f$. If $f\colon (\mathcal Y,S)\to  (\mathcal Z,z)$ is a germ of reflection mapping, a $W$-unfolding of $f$ is the germ at $S\times\{0\}$ of a $W$-unfolding of a representative.
 \end{defn}
 Observe that the triviality condition on the deformation $\widetilde{\mathcal Y}\to \C^r$ ensures that $(\widetilde{\mathcal Y},S\times \{0\})$ and $(\mathcal Y,S)\times (\Delta,0)$ are isomorphic. In particular, every $W$-unfolding of a germ $f$ is, up to $\mathcal A$-equivalence, an unfolding of $f$.
 
\begin{ex}[A family of tetrahedral reflection mappings]\label{exTetrahedralReflectionMappings}
	Let $\mathfrak S_4$ be the group of symmetries of a tetrahedron, as in Example \ref{examS4Group}. Consider 
	the family of reflection mappings
	\[f_t^{\mathfrak S_4}=\omega\vert_{\mathcal Y_t}\colon \mathcal Y_t\to \C^3\quad\text{with}\quad\mathcal Y_t=\{u=t(2v+w)\}.\]
	Equivalently, we may parametrize $\mathcal Y_t$ by $(x,y)\mapsto (t(2x+y),x,y)$ and think of $f_t^{\mathfrak S_4}$ as mappings $\C^2\to \C^3$ of the form
	\[\big(x,y\big)\mapsto \Big(t^2(2x+ y)^2+x^2+y^2,\big(t^2(2x+ y)^2-x^2\big)y,\big(2t^2(2x+ y)^2-y^2\big)\big(2x^2-y^2\big)\Big)\]
	
	As it turns out, $\mathcal Y_0=\{u=0\}=\Fix(1\; 2)$ is one of the reflecting hyperplanes of  $\mathfrak S_4$. Since $\mathfrak S_4$ acts transitively on its reflecting hyperplanes, we have that $W\mathcal Y_0=\mathscr A$, thus \[\im f_0^{\mathfrak S_4}=\omega(\mathcal Y_0)=\omega(\mathscr A),\]
	which means that the image of $f_0^{\mathfrak S_4}$ is precisely the discriminant of the orbit map $\omega$.
\end{ex}

\subsection*{The degree of a reflection mapping}
Before turning our attention to the image and double points space, we show how the degree of a reflection map is encoded by the stabilizers of $\mathcal Y$ in a very simple way. This is essential to our study of the image in the hypersurface case, but it applies to all dimensions.

In this work, when we talk about degree of a mapping $f\colon \mathcal Y\to \mathcal Z$, we mean the number of preimages of a generic point in $f(\mathcal Y)$. For this to make sense, there must be an open dense subset of $f(\mathcal Y)$ on which this number of preimages is contant. This happens if $f(\mathcal Y)$ is irreducible and $f\colon\mathcal Y\to f(\mathcal Y)$ is proper. The degree for finite map-germs is defined by taking a proper representative. For mappings where $f(\mathcal Y)$ fails to be irreducible, a different degree is associated to $f$ on each irreducible component of $f(\mathcal Y)$.


\begin{prop}\label{propDegreeReflectionMapping}
	Let $f$ be a reflection mapping such that $f\colon\mathcal Y\to f(\mathcal Y)$ is proper and $f(\mathcal Y)$ is irreducible. The degree of $f$ is $|W^\mathcal Y:W_\mathcal Y|$.
	\begin{proof}
		The degree of $f$ is the number of preimages of a generic point in $\omega(\mathcal Y)$, that is, the number of points in $Wu\cap \mathcal Y$, for $u$ in a certain open dense subset $\mathcal U_1\subseteq \mathcal Y$. Observe that, given $\sigma\in W$ and $u\in \mathcal Y$, the condition $\sigma u\in \mathcal Y$ may hold even if $\sigma$ does not fix the set $\mathcal Y$, that is, even if $\sigma\notin \mathcal W^{\mathcal Y}$. But, by definition of $\mathcal W^{\mathcal Y}$, this must happen only on a proper closed subset of $\mathcal Y$. Hence, there is an open dense subset $\mathcal U_2\subseteq \mathcal Y$ such that the orbit of $u\in \mathcal U_2$ by the action of $W^\mathcal Y/W_\mathcal Y$ is  $Wu\cap \mathcal Y$. Since  $W^\mathcal Y/W_\mathcal Y$ acts faithfully on $\mathcal Y$, there is an open dense subset $\mathcal U_3$, such that the orbit of a point $u\in \mathcal U_3$ consists of $|W^\mathcal Y:W_\mathcal Y|$ points. Since $ \mathcal U_1\cap  \mathcal U_2\cap  \mathcal U_3\neq \emptyset$, the claim follows.
	\end{proof}
\end{prop}

\begin{ex}\label{exDegree}
	Consider the mapping $f(x,y)=(x^k, y^k, xy)$ of Example \ref{exReflectedGraphs}. One sees that the subgroup of $\Z/k\times\Z/k$ preserving  $\mathcal Y=\{u_3=u_1u_2\}$ is \[(\Z/k\times\Z/k)^\mathcal Y=\langle(1,-1)\rangle\cong \Z/k,\] while no non trivial element of $\Z/k\times\Z/k$ preserves $\mathcal Y$ pointwise. Therefore, the mapping $f$ is generically $k$-to-one. The mapping $f_0^{\mathfrak S_4}$ of Example \ref{exTetrahedralReflectionMappings} has
	\[{\mathfrak{S}_4}^{\mathcal {Y}_0}=\langle (1\,2) ,(3\,4) \rangle\cong \Z/2\times\Z/2\quad \text{\and}\quad {\mathfrak S_4}_{\mathcal {Y}_0}=\langle (1\,2) \rangle\cong \Z/2.\]
	Therefore, $f_0^{\mathfrak S_4}$ parametrizes the discriminant of $\omega$ in a two-to-one way. 
\end{ex}

\begin{cor}\label{corGenericallyOneToOne}
	Let $f$ be a reflection mapping such that $f\colon\mathcal Y\to f(\mathcal Y)$ is proper. Then, $f$ is generically one-to-one if and only if, for all $\sigma\in W\setminus\{1\}$, $\dim((\mathcal Y\cap\sigma\mathcal Y)\setminus \Fix\sigma)\leq n-1$.
	\begin{proof}We may assume $\omega(\mathcal Y)$ to be irreducible, for neither the generically one-to-one property nor the condition $\dim((\mathcal Y\cap\sigma\mathcal Y)\setminus \Fix\sigma)\leq n-1$ will be affected if we consider the irreducible components of $\omega(\mathcal Y)$ separately.   The condition $\dim((\mathcal Y\cap\sigma\mathcal Y)\setminus \Fix\sigma)\leq n-1$ is equivalent to the statement that any $\sigma\in W$ must satisfy either $\mathcal Y\neq \sigma\mathcal Y$ or $\mathcal Y\subseteq\Fix\sigma$, that is, that any $\sigma\in W^\mathcal Y$ must be contained in $ W_\mathcal Y$.
	\end{proof}
\end{cor}

In order to prove our formulas for $\im f$, we need the next result. Since the proof is slightly involved and the statement is quite clear, we postpone the proof until Appendix \ref{secDelayedProofs}.
\begin{lem}[Generically One-To-One Unfolding]\label{lemGen1To1Unfolding}
	Any multi-germ of reflection mapping admits a one-parameter $W$-unfolding, given by $\widetilde{\mathcal Y}\subseteq \mathcal V\times \C$, such that $\dim(\widetilde{\mathcal Y}\cap\sigma\widetilde{\mathcal Y})< \dim\widetilde{\mathcal Y}$, for all $\sigma\in W\setminus\{1\}$.
	In particular, this unfolding is generically one-to-one.
\end{lem}

\subsection*{Fitting ideals}

It is well known that the image of a finite holomorphic mapping is an analytic set (indeed it is enough for the mapping to be proper). However, given an unfolding $F=(f_t,t)$, the ideal of the image of $f_0$ need not be the same as the result of computing the ideal of the image of $F$ and replacing $t=0$. This is a problem for the study of singular mappings, where deformations are regarded as an essential part of the theory. Luckily enough, there is a solution consisting on declaring the image of a finite mapping $f\colon \mathcal Y\to \mathcal Z$ to be, rather than just a set, a complex space
\[\im f=V(\mathcal F_0(f_*\mathcal O_{\mathcal Y})),\] where $\mathcal F_0(f_*\mathcal O_{\mathcal Y})$ stands for the $0$th Fitting ideal sheaf of the pushforward module. This sometimes gives $\im f$ a non-reduced analytic structure, but this is the price we pay in order for the analytic structure to behave well under deformations. For map-germs, one uses the $0$th Fitting ideal, written $F_0(f_*\mathcal O_{\mathcal Y})$. For information about Fitting ideals, we refer to \cite{MondPellikaanFittingIdeals}. We only include here the results we need.
	
	\begin{prop}\label{propFittingUnfolding} If  $F=(f_t,t)\colon \mathcal Y\times \Delta\to  \mathcal Z\times \Delta$ is a finite mapping  and $\epsilon\in \Delta$ then $V(\mathcal F_0(f_{\epsilon*}\mathcal O_{\mathcal Y}))=V(\mathcal F_0(F_*\mathcal O_{\mathcal Y}))\cap \{t=\epsilon\}$.
	\end{prop}

	\begin{prop}\label{propFittingIdealsMultiplicity}
Let $f\colon X\to \C^{n+1}$ be a finite mapping defined on a reduced $n$-dimensional Cohen-Macaulay space, and assume the irreducible decomposition of $X$ to be $X_1,\dots,X_m$. Let $g_i$ be a generator of the ideal of $f(X_i)$ in $\C^{n+1}$, and let $d_i$ be the degree of $f$ restricted to $X_i$. Then $F_0(f_*\mathcal O_X)$ is a principal ideal, generated by 
\[g= \prod_{i=1}^m g_i^{d_i}.\]
	\end{prop}
	
	The first result follows from \cite[Lemma 1.2]{MondPellikaanFittingIdeals}, the second is  \cite[Proposition 3.2]{MondPellikaanFittingIdeals}. The next one is a slight modification, tailored to our needs:
	
	\begin{prop}\label{propFittingIdealsDecomposition}
Let $f\colon X\to \C^{n+1}$ be a finite map-germ defined on an $n$-dimensional Cohen-Macaulay space. Let $X_1,\dots,X_r$ be $n$-dimensional Cohen-Macaulay complex spaces, forming a set theoretical decomposition $X=\cup_{i=1}^rX_i$,
where $X_i$ and $X_j$ have no common components if $i\neq j$.
Assume each $X_i$ to be isomorphic to $X$ on an open dense subset of $X_i$. Then,
\[F_0(f_*\mathcal O_X)=\prod_{i=1}^r F_0((f\vert_{X_i})_* \mathcal O_{X_i}).\]
\begin{proof}
Let $\tilde X$ be the disjoint union of the spaces $X_i$ and take the obvious mapping $\tilde f\colon \tilde X\to \C^{n+1}$. Since being Cohen-Macaulay is a local property satisfied by each of the $X_i$, the space $\tilde X$ is Cohen-Macaulay. Since \[\tilde f_*\mathcal O_{\tilde X}=\bigoplus_{i=1}^r (f\vert_{X_i})_* \mathcal O_{X_i},\] it follows from the construction of Fitting ideals that $F_0(\tilde f_*\mathcal O_{\tilde X})=\Pi_{i=1}^r F_0((f\vert_{X_i})_* \mathcal O_{X_i})$. Since $f$ and $\tilde f$ are the same map on an open dense subset of their targets, it follows that 
the stalks of $F_0(\tilde f_*\mathcal O_{\tilde X})$ and $F_0(f_*\mathcal O_{X})$ are the same on that open dense subset. But $F_0(\tilde f_*\mathcal O_{\tilde X})$ and $F_0(f_*\mathcal O_{X})$ are principal ideals, hence they must agree.\end{proof}

	\end{prop}

\section{The image of a reflection mapping $\mathcal Y^n\to \C^{n+1}$}\label{secImage}

	In this section we show explicit formulas for the image of a reflection map $\mathcal Y^n\to \C^{n+1}$. For any finite map-germ $f\colon  (\mathcal Y^{n},0)\to (\mathcal Z^{n+1},0)$ between complex manifolds, the ideal $F_0(f_*\mathcal O_{\mathcal Y})$ is principal (this holds, more generally, whenever $\mathcal Y$ is  an $n$-dimensional Cohen Macaulay space, and it follows from \cite[Section 2.2]{MondPellikaanFittingIdeals}). Hence, letting $g$ be a generator of $F_0(f_*\mathcal O_{\mathcal Y})$, the image of $f$ is
\[\im f=V(g).\]

Putting together Proposition \ref{propFittingIdealsMultiplicity} and Corollary \ref{corGenericallyOneToOne}, one obtains what follows:

	\begin{prop}\label{corReducedImage}
For any reflection mapping $\mathcal Y^n\to \C^{n+1}$, the space $\im f$ is reduced if and only if $f$ is generically one-to-one, if and only if $\dim((\mathcal Y\cap\sigma\mathcal Y)\setminus \Fix\sigma)\leq n-1$,  for all $\sigma\in W\setminus\{1\}$, .
	\end{prop}

	\begin{thm}\label{thmImageEquation}
For any reflection mapping $\mathcal Y^n\to \mathcal \C^{n+1}$, $\im f$ is the zero locus of
\[g=\prod_{\sigma\in W}(\sigma L)\circ s,\]
where $s\colon \C^{n+1}\to\mathcal V$ is a section of $\omega$. Equivalently, $g(X)=\Pi_{u\in\omega^{-1}(X)}L(u)^{|W^u|}$.
\begin{proof}
It is immediate that $g$ vanishes precisely at $\im f$, hence we must check that $g$ is holomorphic and that the analytic structure of $V(g)$ is that of $V(F_0(f_*\mathcal O_{\mathcal Y}))$. This may be verified locally in the target and, for simplicity, we will do it just at the origin. 

Before comparing both structures, we show that $g$ is holomorphic and study the space $V(g)$. The ring homomorphism $\omega^*\colon \mathcal O_{n+1}\to\mathcal O_{\mathcal V,0}$ may be identified with the inclusion $\mathcal O_{\mathcal V,0}^W\hookrightarrow \mathcal O_{\mathcal V,0}$, where $\mathcal O_{\mathcal V,0}^W$ stands for the subring of $W$-invariant germs. This identifies the ideal $\omega_*^{-1}(\langle\Pi_{\sigma \in W}\sigma L)\rangle$ in $\mathcal O_{n+1}$ with $\langle \Pi_{\sigma \in W}\sigma L\rangle \cap \mathcal O_{\mathcal V,0}^W$. From the fact that $\langle \Pi_{\sigma \in W}\sigma L\rangle$ is a principal ideal generated by a $W$-invariant function, it follows that this ideal is precisely the ideal generated by $ \Pi_{\sigma \in W}\sigma L$ in $\mathcal O_{\mathcal V,0}^W$. In particular, $\omega_*^{-1}(\Pi_{\sigma \in W}\sigma L)$ is a principal ideal generated by a holomorphic function $\tilde g$, such that  $\Pi_{\sigma\in W} \sigma L=\tilde g\circ \omega$. Since $\omega$ is surjective, the function $\tilde g$ is uniquely determined, not only among the holomorphic functions, but among all functions $\C^{n+1}\to \C$.  Since the function $g(X)=\Pi_{u\in\omega^{-1}(X)}L(u)^{|W^u|}$
satisfies the desired equality, we conclude $g=\tilde g$ and, in particular, $g$ is holomorphic. Moreover, we have identified the coordinate ring of $V(g)$ with the quotient of $\mathcal O_{\mathcal V}^W$ by  $\langle \Pi_{\sigma \in W}\sigma L\rangle \cap \mathcal O_{\mathcal V,0}^W$. Note that this ring is a subring of $\mathcal O_\mathcal V/\langle \Pi_{\sigma \in W}\sigma L\rangle$.

Now we show $V(g)\cong V(F_0(f_*\mathcal O_\mathcal Y))$. Obviously, the explicit formula given for $g$ behaves well under $W$-unfoldings in the same way $F_0(f_*\mathcal O_\mathcal Y)$ does for any unfolding (see Proposition \ref{propFittingUnfolding}). Consequently, it suffices to show the isomorphism for an unfolding of $f$. Hence, in view of the Generically One-to-one Unfolding Lemma \ref{lemGen1To1Unfolding}, we may assume $f$ to satisfy $\dim(\mathcal Y\cap\sigma\mathcal Y)\leq n-1$. Since  $V(g)$ and $V(F_0(f_*\mathcal O_\mathcal Y))$ are the same at the set theoretical level, it suffices for both spaces to be reduced  for them to be equal. The space $V(F_0(f_*\mathcal O_\mathcal Y))$ is reduced by Proposition \ref{corReducedImage}. At the same time the condition $\dim(\mathcal Y\cap\sigma\mathcal Y)\leq n-1$ is equivalent to the condition that $V( \Pi_{\sigma \in W}\sigma L)$ is reduced, which forces $V(g)$ to be reduced because, as we mentioned before, the stalk of $\mathcal O_{V(g)}$ is a subrings of stalks of $\mathcal O_\mathcal V/\langle \Pi_{\sigma \in W}\sigma L\rangle$.
\end{proof}
	\end{thm}
	
	\begin{rem}The same formula applies to compute images of \emph{singular} hypersurfaces. To be precise,
if $\mathcal Y=V(L)\subseteq \mathcal V$ is a singular hypersurface (instead of a complex manifold), then $g=\Pi_{\sigma\in W}(\sigma L)\circ s$ computes a generator of the ideal $F_0((\omega_{\vert_\mathcal Y})_*\mathcal O_{\mathcal Y})$, which defines the image of $\mathcal Y$ by $\omega$. Briefly speaking, the formula holds because it is correct on the dense open subset where $\mathcal Y$ is smooth, but this forces it to be correct everywhere. If $\mathcal Y$ is non reduced, one shows the claim by taking a $W$-unfolding with reduced generic fiber.
	\end{rem}

\begin{ex}Consider the group $\Z/d_1\times\dots\times \Z/d_{n+1}$ acting on $\C^{n+1}$. Then,
\[g=\prod_{\substack{0\leq a_i< d_i}}L(\xi_1^{a_1}X^{\frac{1}{d_1}},\dots, \xi_{n+1}^{a_{n+1}}X^{\frac{1}{d_{n+1}}} )\]
One does not need to worry about the definition of the complex square root, because any choice of a section of $\omega$ will give the same final result. 

Take, for example the fold mappings  $(x,y)\mapsto (x,y^2,yP(x,y^2))$ of Example \ref{exReflectedGraphs}. The ideal of $\mathcal Y$ is generated by $L(u,v,w)=w-vP(u,v^2)$, hence the image is the zero locus of
\[g=(Z-\sqrt{Y}P(X,Y))(Z+\sqrt{Y}P(X,Y))=Z^2-YP^2(X,Y).\]
For double folds  $(x,y)\mapsto (x^2,y^2,xP_1(x^2,y^2)+yP_2(x^2,y^2)+xyP_3(x^2,y^2))$ (also Example \ref{exReflectedGraphs}) one has $L(u,v,w)=w-uP_1(u^2,v^2)-vP_2(u^2,v^2)-uv P_3(u^2,v^2)$ and,  writing $P_i=P_i(X,Y)$, one gets
\begin{align*}
g=&(Z-\sqrt{X}P_1-\sqrt{Y}P_2-\sqrt{XY}P_3)\cdot(Z+\sqrt{X}P_1-\sqrt{Y}P_2+\sqrt{XY}P_3)\cdot\\
	&(Z-\sqrt{X}P_1+\sqrt{Y}P_2+\sqrt{XY}P_3)\cdot(Z+\sqrt{X}P_1+\sqrt{Y}P_2-\sqrt{XY}P_3)\\
	=&(Z^2-XP_1^2-YP_2^2-XYP_3^2)^2-4XY(2P_1P_2P_3+P_1^2P_2^2+XP_1^2P_3^2+YP_2^2P_3^2).
\end{align*}

 \end{ex}

%
%

\begin{ex}[Image of $f^{ D_8 }$]\label{d8image}
Consider the  $D_8$-reflected graph of Example \ref{exReflectedGraphs},
\[f^{ D_8 }_1\colon (x,y)\mapsto (x^2+y^2, x^2y^2, 2x+y).\]
A generator of the ideal of $\mathcal Y$ is $L(u,v,w)=w-2u-v$ and a section of  $\omega$ is, for example,
\[s(X,Y,Z)=\left(\sqrt{\frac{2Y}{X-\sqrt{X^2-4Y}}}, \sqrt{\frac{X-\sqrt{X^2-4Y}}{2}},Z\right).\]

 Multiplying the elements in the orbit of $L=w-2u-v$ gives the function
\begin{align*}
\prod_{\sigma\in W}\sigma L=&(w-2u-v)(w+2u-v)(w-2u-v)(w-2u+v)\\&(w+2v+u)(w-2v+u)(w+2u+v)(w+2v-u).
\end{align*} 
Taking the composition  with $s$, the square roots vanish and we obtain the expression
\[g=16X^4-200X^2Y+625Y^2-40X^3Z^2+70XYZ^2+33X^2Z^4-14 YZ^4-10 XZ^6+Z^8.\]

The same method gives the equation of the image of 
\[f^{ D_8 }_2\colon (x,y)\mapsto (x^2+y^2,x^2y^2,2x^2+3xy-y^2+2x^3+8x^2y-2xy^2-2y^3),\] but it is too big to be written here. Both surfaces are depicted in Figure \ref{imaged82}.

\end{ex}

\begin{rem}\label{310}
Let $f=(\omega,H)$ be a $W$-reflected graph, and let $d_1,\dots,d_n$ be the degrees of the coordinates functions of $\omega$. If $H$ is homogeneous of degree $d$, then  $f$ is obviously homogeneous, with degrees $d_1,\dots,d_n,d$. As a consequence of Theorem \ref{thmImageEquation}, the function $g$ defining $\im f$ is quasi-homogeneous with weights
$d_1,\dots,d_n,d$ and degree $d|W|$.
\end{rem}

\subsection*{Computing the image  via elimination of variables}
Computing explicit expressions of sections of orbit mappings can be hard and, even if we manage to find them, the formula $\Pi_{\sigma\in W}\,\sigma L\circ s$ becomes too tedious to compute by hand when we take bigger reflection groups. If we try to implement the expression in some computer software, we face the fact that the coordinate functions of sections are non-polynomial, which is a problem when using commutative algebra software such as {\sc Singular}. Here we introduce two alternative ways of computing the ideal $F_0(f_*\mathcal O_\mathcal Y)$ defining the image, without resorting to sections of $\omega$. As a trade-off, this methods use elimination of variables, which is less explicit, but it is a task computers are happy to do, or at least to try to.

In the case at hand, elimination of variables means what follows: Consider the ring $\C[u,X]$ of polynomials on the variables $u_1,\dots,u_p$ and $X_1,\dots X_p$ (the same works for, say, the ring of germs of holomorphic functions at the origin). Contained in $\C[u,X]$, there is the ideal
\[J=\langle\prod_{\sigma\in W}\sigma L(u)\rangle+\langle X_1-\omega_1(u),\dots,X_p-\omega_p(u)\rangle,\]
and the subring $\C[X]$ of polynomials in $X_1,\dots X_p$. Then, $J\cap \C[X]$ is an ideal in $\C[X]$, said to be obtained by \emph{eliminating} the $u$ variables in $J$. Alternatively, consider the ring homomorphism
\[\omega^*\colon \C[X]\to \C[u],\]
given by $H\mapsto H\circ \omega$. Then, the elimination can be expressed as \[J\cap \C[X]=(\omega^*)^{-1}\langle\prod_{\sigma\in W}\sigma L(u)\rangle.\]

\begin{prop}
For any reflection mapping $f\colon \mathcal Y^n\to\C^{n+1}$, $\im f$ is the zero locus of
\[F_0(f_*\mathcal O_\mathcal Y)=(\omega^*)^{-1}\langle \prod_{\sigma\in W}\sigma L \rangle=\left((\omega^*)^{-1}\langle  L \rangle\right)^{|W^\mathcal Y:W_\mathcal Y|}.\]
\begin{proof}The first equality was shown in the  proof of Theorem \ref{thmImageEquation}. For the second equality, let 
\[K=\langle L(u)\rangle+\langle X_1-\omega_1(u),\dots,X_p-\omega_p(u)\rangle,\]
Since $V(K)\cong V(L)=\mathcal Y$, it follows that $K$ is radical, hence $K\cap \C[X]$ is radical. This means that $K\cap \C[X]$ is the (radical) ideal of the image of $f$. At the same time, we know by Proposition \ref{propDegreeReflectionMapping} that the degree of $f$ is $|W^\mathcal Y:W_\mathcal Y|$. Putting both things together, the equality $F_0(f_*\mathcal O_\mathcal Y)=(K\cap \C[X])^{|W^\mathcal Y:W_\mathcal Y|}$ is a direct application of Proposition \ref{propFittingIdealsMultiplicity}.
\end{proof}
\end{prop}

In practice, our experience shows that computing $(\omega^*)^{-1}\langle \Pi_{\sigma\in W}\sigma L \rangle$ is faster than computing $(\omega^*)^{-1}\langle  L \rangle$, probably due to the fact that $\Pi_{\sigma\in W}\sigma L$ is a $W$-invariant function.
\begin{ex}
Consider the family of $\mathfrak{S_4}$-reflection mappings of Example \ref{exTetrahedralReflectionMappings}. Sections of the orbit mapping of $\mathfrak{S_4}$ have horrendous expressions, so computing $\im f$ by means of the expression $\Pi_{\sigma\in W}\,\sigma L\circ s$ of Theorem \ref{thmImageEquation} is not convenient. In contrast, {\sc Singular} computes $(\omega^*)^{-1}\langle \Pi_{\sigma\in W}\sigma L \rangle$ in no time. The equation of the image of the unfolding $F=(f_t,t)$ is too long to write here but, for $t=0$ and $t=1$, we obtain the following equations:
\[\im f^{\mathfrak{S_4}}_0=V(g),\quad g=(2 x^3 y^2+ x^4 z-27 y^4-18  x y^2 z-2  x^2 z^2+z^3)^2,
\]
\[
\begin{split}
\im f^{\mathfrak{S_4}}_1=V(g),\quad g=614656x^{12}-174822592x^9y^2-16020256x^{10}z+10356692964x^6y^4\\
+800288220x^7y^2z+153738321x^8z^2-198333009364x^3y^6-33901243950x^4y^4z\\
-662345364x^5y^2 z^2-685828516 x^6 z^3+1202174306137 y^8+372758486548 x y^6 z\\
+7876328208 x^2 y^4 z^2-1163406956 x^3 y^2 z^3+1546928326 x^4 z^4+40000919994 y^4 z^3\\
+1284226020 x y^2 z^4-1713759300 x^2 z^5+741200625 z^6.
\end{split}
\]
The hypersurces $\im f^{\mathfrak{S_4}}_0$ and $\im f^{\mathfrak{S_4}}_1$ are depicted in Figure \ref{figBishopAndBat}. Recall that $\im f^{\mathfrak{S_4}}_0$ is the discriminant of $\omega$, with non-reduced structure, and that the exponent $2$ in the equation reflects the fact that $f^{\mathfrak{S_4}}_0$ has degree two (see Example \ref{exDegree}).
 \begin{figure}[h]
\begin{center}
\includegraphics[scale=.6]{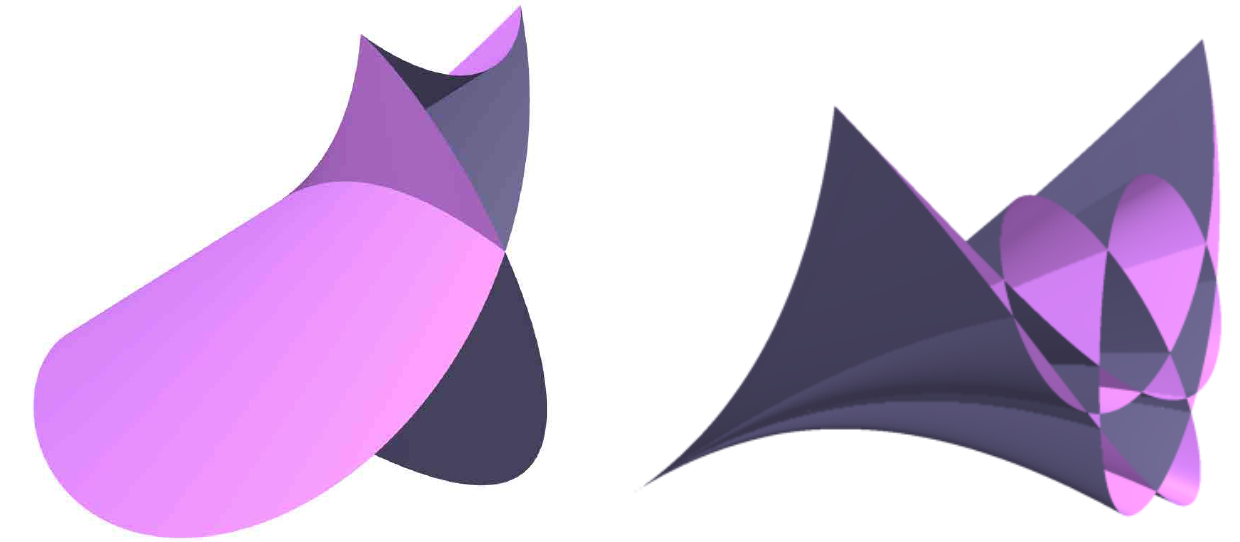}
\end{center}
	\caption{The images of $f_0^{\mathfrak S_4}$ and 
$f_1^{\mathfrak S_4}$.}
	\label{figBishopAndBat}
\end{figure}
\end{ex}

\begin{rem}
The calculations shown before have been carried out by means of a  {\sc Singular} library for reflection mappings, being developed by the authors of this paper. We are happy to send the latest version to anyone interested.
\end{rem}

\section{Decomposition of the  double point spaces $K_2(f), D^2(f)$ and $D(f)$.}
Given a germ of singular mapping, studying only the equations of the image somehow means forgetting that these are parametrizable singularities, that is, that these singularities are obtained by glueing one disk by means of a holomorphic mapping (or several disks, in the multi-germ case). Instead, it usual to study singular mappings by means of their multiple-point spaces, which are spaces designed to encode how points in the source glue together to form the image of the map. 

 There are double point and triple point spaces, as well as higher multiplicity ones, but here we restrict ourselves to double points. Even then, there are several different double point spaces one can look at. Our original interest is to prove an explicit formula  for the double point space $D(f)\subseteq \mathcal Y$ in the case where $f$ is a reflection mapping given by a hypersurface $\mathcal Y\subseteq \mathcal V$ (see Theorem \ref{thmAFormula}. The definitions of $D(f)$ and other double point spaces are given later on this section). However, as the logic structure of this and the following section reveal, this is best done by studying first the double point space $D^2(f)\subseteq \mathcal Y\times\mathcal Y$. This in turn requires looking at a more abstract double point space $K^2(f)$, which enjoys a key functorial property (Proposition \ref{propFunctorialityOfK2}). The results before Section \ref{secAFormula} do not require $\mathcal Y$ to be a hypersurface and are given for all codimensions.

\subsection*{Decomposition of the double point space $K_2(f)$}
Here we study Kleiman's double points for reflection mappings. We skip many details, for which we refer to \cite{NunoPenafortIterated, MR309129}.
For this description of the double point space $K^2(f)$ for   mappings (not necessarily reflection mappings) let $\mathcal Y$ be an $n$-dimensional complex manifold and assume that it admits a global coordinate system (This is for the sake of simplicity, a standard glueing process shows the results in this section to be valid in general). We write the blowup of the product of two copies of $\mathcal Y$ along the diagonal as
\[B_2(\mathcal Y)=\operatorname{Bl}_{\Delta\mathcal Y}(\mathcal Y\times \mathcal Y)=\{(u,u',v)\in \mathcal Y\times \mathcal Y\times \mathbb P^{n-1}\mid v\wedge(u'-u)=0\},\]
where $v\wedge(u'-u)=0$ is a shortcut to indicate the vanishing of all $2\times 2$ minors of the matrix
\[
\left(\begin{array}{ccc}
u'_1-u_1 & \dots & u'_n-u'_n \\
v_1 & \dots & v_n
\end{array}\right).
\]
We write the exceptional divisor as $E=\Delta\mathcal Y\times  \mathbb P^{n-1}$.

Now let $f\colon \mathcal Y\to \mathcal Z$ be a holomorphic mapping between manifolds (not necessarily a reflection mapping) of dimensions $n$ and $p$, both admiting a global coordinate system.   Think of $u'-u$ and $f(u')-f(u)$ as a vectors with entries in $\mathcal O_{\mathcal Y\times \mathcal Y}$, of sizes $n$ and $p$, respectively. By Hilbert Nullstellensatz, there exist a $p\times n$ matrix $\alpha_f$, with entries in $\mathcal O_{\mathcal Y\times \mathcal Y}$, such that
\[f(u')-f(u)=\alpha_f\cdot (u'-u)\]
Then, Kleiman's double point space is
\[K_2(f)=\{(u,u',v)\in B_2(\mathcal Y)\mid \alpha_f(u,u')\cdot v=0\}.\]
The proof that $K_2(f)$ does not depend on $\alpha$ is in \cite[Proposition 3.1]{MondSomeRemarks}. Away from the exceptional divisor, $K_2(f)$  is just the fibered product \[K_2(f)\setminus E \cong (\mathcal Y\times_{\mathcal Z}\mathcal Y)\setminus \Delta\mathcal Y,\] via the blowup map.
On the exceptional divisor, $K_2(f)$ keeps track of the kernel of the differential of $f$, as follows:
\[K_2(f)\cap E=\{(u,u,v)\in B_2(\mathcal Y)\mid v\in \ker \operatorname{d}\!f_u\}.\]

\begin{prop}\label{propK2CompleteIntersection}\cite[Corollary 5.6]{NunoPenafortIterated}
The dimension of $K_2(f)$ is at least $2n-p$. If $K_2(f)$ has dimension $2n-p$, then it is locally a complete intersection.
\end{prop}

The double point spaces $K_2$ enjoy the following key functorial property, very much related to the construction of reflection mappings. First, any embedding of complex manifolds, $\mathcal Y\hookrightarrow \mathcal V$, induces an embedding $B^2(\mathcal Y)\hookrightarrow B^2(\mathcal V)$. 

\begin{prop}\cite[Theorem 2.13]{NunoPenafortIterated}\label{propFunctorialityOfK2}
 Given a mapping $F\colon \mathcal V\to \mathcal Z$ between complex manifolds, the double point space of the restriction $f=F_{\vert_\mathcal Y}\colon \mathcal Y\to \mathcal Z$ is
\[K_2(f)=K_2(F)\cap B_2(\mathcal Y).\]
\end{prop}
In particular, the double points of reflection mappings are slices of the double points of orbit mappings. Conveniently enough, we understand $K_2(\omega)$ quite well. Fixed $\sigma \in W\setminus\{\id\}$, let $r=\dim \Fix \sigma$ and choose linear mappings $\ell_\sigma\colon \mathcal V\to \C^{p-r}$ and $\ell_\sigma^\bot\colon \mathcal V\to \C^{r}$, so that
\[\Fix \sigma=V(\ell_\sigma)\qquad\text{and}\qquad(\Fix \sigma)^\bot=V(\ell_\sigma^\bot).\]
The following description of $K^2(\omega)$ is found (with different notation) in \cite[Theorem 7.6]{Penafort-Sanchis:2016a}:

\begin{prop}  The double point space $K_2(\omega)$ is a $p$-dimensional reduced locally complete intersection, with irreducible decomposition
\[K_2(\omega)=\bigcup_{\sigma\in W\setminus\{1\}}K_2^\sigma,\]
where $K_2^\sigma$ is the blow-up of  $\{(u,\sigma u)\mid u\in \mathcal V\}$ along $\Delta \Fix \sigma$, embedded in $B_2(\mathcal V)$ as
\[K_2^\sigma=\{(u,\sigma u,v)\in B^2(\mathcal V)\mid  \ell_\sigma^\bot(v)=0\}.\]
\end{prop}

In view of Proposition \ref{propFunctorialityOfK2}, the space $K_2$ of reflection mappings must inherit a $W\setminus\{1\}$-indexed decomposition:
\begin{defn} For each $\sigma\in W\setminus\{1\}$, we let $K_2^\sigma(f)=B_2(\mathcal Y)\cap K_2^\sigma.$
\end{defn}
As we just mentioned, this gives  the set-theoretical decomposition
\[K_2(f)=\bigcup_{\sigma\in W\setminus\{1\}}K_2^\sigma(f).\]
Observe that neither $K_2(f)$ nor $K_2^\sigma(f)$ need to be reduced, so there is no clear way to upgrade the previous decomposition into a union of complex spaces, because (to the best of the authors' knowledge) the meaning of a union means in the category of complex spaces is unclear. In any case, since the order of $W$ is finite,  $K_2(\omega)$ and $K_2^\sigma$ are locally equal topological spaces at points $z\in K_2^\sigma$ not contained in any other $K_2^\tau,\tau \neq \sigma$. Moreover, since both spaces are reduced, $K_2(\omega)$ and $K_2^\sigma$ are locally isomorphic complex spaces at such points. This means that the comparison of analytic structures in the above equality is only problematic  at points where different branches $K_2^\sigma(f)$  and  $K_2^\tau(f)$ meet. To be precise, we have observed what follows:

	\begin{lem}\label{lemK2K2SigmaAwayFromK2Tau}
Let $z\in K_2(f)$ be a point contained in $K_2^\sigma(f)$. If $z$ is not contained in $K_2^\tau(f)$, for any $\tau\neq \sigma$, then $K_2(f)$ and $K_2^\sigma(f)$ are locally isomorphic at $z$.
	\end{lem}
	
	Now we want to refine our description of $K_2^\sigma(f)$.
Again, we may think of  $L$ and $\sigma L$ as vectors with entries in $\mathcal O_{\mathcal V}$. Then, there is a $(p-n)\times (p-r)$ matrix $\alpha_L^\sigma$, with entries in $\mathcal O_{\mathcal V}$, such that 
\[\sigma^{-1}L-L=\alpha_L^\sigma \cdot \ell_\sigma.\]

Also, consider the involution $\iota\colon \mathcal V\times \mathcal V\to \mathcal V\times \mathcal V$, given by $(u,u')\mapsto (u',u)$. This mapping lifts to an involution $\tilde \iota \colon B^2(\mathcal V)\to B^2(\mathcal V)$  taking $K_2(f)$ onto itself.

	\begin{prop}\label{propK2AnalyticStructure}
With the above notations, the branch $K^\sigma_2(f)$ is the complex space
\[K_2^\sigma(f)=\{(u,\sigma u,v)\in K^\sigma_2 \mid  L(u)=0, \alpha_L^\sigma(u)\cdot \ell_\sigma(v)=0\}\]
and it satisfies the following properties:
\begin{enumerate}
\item \label{propK2AnalyticStructureItem1} The dimension of $K_2^\sigma(f)$ is at least $2n-p$. If $ K_2^\sigma(f)$ has dimension $2n-p$, then it is locally a complete intersection.
\item \label{propK2AnalyticStructureItem2}The blowup mapping $B_2(\mathcal Y)\to \mathcal Y\times \mathcal Y$ takes $K_2^\sigma(f)\setminus E$ isomorphically to the space
\[\{(u,\sigma u)\mid u\in (\mathcal Y\cap \sigma^{-1}\mathcal Y)\setminus\Fix\sigma\}.\]
\item \label{propK2AnalyticStructureItem3}The involution $\tilde \iota$ takes $K_2^\sigma(f)$ isomorphically to $K_2^{(\sigma^{-1})}(f)$.
\end{enumerate}
\begin{proof}
The space in question is defined as $K_2^\sigma(f)=B^2(\mathcal Y)\cap K_2^\sigma$. Consider the $(p-n)\times p$ matrix $\alpha_L$, with entries in $\mathcal O_{\mathcal V\times \mathcal V}$, satisfying $L(u')-L(u)=\alpha_L(u,u')\cdot(u'-u)$. This allows to describe $B_2(\mathcal Y)$ inside $B_2(\mathcal V)$ as 
\[B_2(\mathcal Y)=\{(u,u',v)\in B_2(\mathcal V)\mid L(u)=0,\alpha_L(u,u')\cdot v=0\}.\] To see this, it suffices to check that the space on the right hand side of the equality is the strict transform of $\mathcal Y\times \mathcal Y$ in $B_2(\mathcal V)$. This in turn is true because this space has dimension $2n$, making it a complete intersection by counting equations, and  it is isomorphic to $\mathcal Y\times \mathcal Y$ away from $E$, while the equations prevent it from having irreducible components contained in $E$. Now, to finish the description of $K_2^\sigma(f)$ in the statement, it suffices to show that, for points $(u,\sigma u,v)\in K_2^\sigma$, the conditions $\alpha_L(u,\sigma u)\cdot v=0$ and $\alpha_L^\sigma(u)\cdot \ell_\sigma(v)=0$ are equivalent. For simplicity, we may identify $\mathcal V=\Fix \sigma^\bot\oplus\Fix \sigma$, so that any vector is of the form $u=(u_1,u_2)$, with  $\ell_\sigma(u)=u_1$ and $\ell_\sigma^\bot(u)=u_2$. Any point $(u,\sigma u,v)\in K_2^\sigma$ satisfies $\ell_\sigma^\bot(v)=0$, which means that the equation $\alpha_L(u,\sigma u)\cdot v=0$ has the form
\[
  \left(\begin{array}{@{}c|c@{}}
 A(u,\sigma u)& B(u,\sigma u)   \end{array}\right)
  \left(\begin{array}{c}
    v_1 \\\hline
    0  \end{array}\right)=0.
\]
Taking $\alpha_L^\sigma(u)=A(u,\sigma u)$, the equation $\alpha_L^\sigma(u)\cdot \ell_\sigma(v)=0$ is equivalent to the one above. 

Item (\ref{propK2AnalyticStructureItem1}) is an easy equation counting: $K_2^\sigma$ is a $p$-dimensional manifold, on which we impose $2(p-n)$ equations. To see item (\ref{propK2AnalyticStructureItem2}), observe that $(u,\sigma u)$ is contained in $E$ if and only if $u\in Fix \sigma$, and that away from $E$, the equations $\alpha_L^\sigma(u)\cdot \ell_\sigma(v)=0$ and $L(u')-L(u)=0$ are equivalent. For item (\ref{propK2AnalyticStructureItem3}), we already know that $\tilde \iota$ is an automorphism on $K_2(f)$ and, since it has the form $(u,\sigma u,v)\mapsto(\sigma u, u, v)$, it must take $K_2^\sigma(f)$ to $K_2^{(\sigma^{-1})}(f)$. This must happen isomorphically, since otherwise $\tilde \iota$ would fail to be an automorphism.
\end{proof}
	\end{prop}

\begin{ex}\label{k2c3c5}Consider the  map-germ
		\[f^{(d_1,\dots,d_5)}\colon(x,y,z)\mapsto (x^{d_1}, y^{d_2},z^{d_3},(x+y+z)^{d_4},(x-y+2z)^{d_5}),\]
 from Example \ref{exRefMapsUnbounded}. The equations of $\mathcal Y$ are $L=(u_1+u_2+u_3-u_4,u_1-u_2+2u_3-u_5)=0$.
		To illustrate the calculation of $K_2^\sigma(f)$, take, for instance, the elements of the form $\sigma_i=(a_1,\dots,a_i,0,\dots,0)$, $i=1,\dots,5$, with $0<a_i< d_i$ (see Example \ref{exProductOfCyclicGroups} for notation).  The pointwise fixed spaces are given by the vanishing of $\ell_{\sigma_i}(u)=(u_1,\dots,u_i)$. From the equations $\sigma^{-1}L-L=\alpha_L^\sigma \cdot \ell_\sigma$ one concludes that
		$\alpha_L^{\sigma_i}$ is the matrix defined by the first $i$ columns of
		\[\alpha_L^{\sigma_5}=
		\left(\begin{array}{ccccc}\xi_1^{a_1}-1&\xi_2^{a_2}-1&\xi_3^{a_3}-1&-1(\xi_4^{a_4}-1)&0 \\
			\xi_1^{a_1}-1&-1(\xi_2^{a_2}-1)&2(\xi_3^{a_3}-1)&0&-1(\xi_5^{a_5}-1)
		\end{array}\right)
		\]
		Since, for $i=1,2$, $\alpha_L^{\sigma_i}$ is injective, the equation  $\alpha_L^{\sigma_i} \cdot \ell_{\sigma_i}(v)=0$ forces $\ell_{\sigma_i}(v)=0$. This, added to the condition $\ell_{\sigma_i}^\bot(v)=0$ for points in $(u,\sigma u,v)\in K_2^{\sigma_i}$, forces $v=0$ and, since a point $v\in \mathbb P^{p-1}$ cannot have all its coordinates equal to zero, we conclude $K_2^{\sigma_i}(f)=\emptyset$.
$K_2^{\sigma_3}(f)$ is a line, with no component on the exceptional divisor, and it is a complete intersection. $K_2^{\sigma_4}(f)$ is a line, with an irreducible component on the exceptional divisor.
			$K_2^{\sigma_5}(f)$ is the union of a line and a 2-dimensional component lying on the exceptional divisor. This branch is not a complete intersection (indeed it is not even Cohen Macaulay).
\end{ex}

\subsection*{Decomposition of the double point space $D^2(f)$}  Now we study a different double point space. As a set, $D^2(f)\subseteq \mathcal V\times\mathcal V$ is just the projection of $K_2(f)$. As a complex space, this space was introduced in \cite{MondSomeRemarks}. The results for which we give no proof can be found in \cite{Nuno-Ballesteros2015On-multiple-poi}. Let $f\colon \mathcal Y\to \mathcal Z$ be a holomorphic mapping between complex manifolds, as in the previous section. We define
\[D^2(f)=V((f\times f)^*(I_\Delta)+I_n(\alpha_f)),\]
where $I_\Delta$ stands for ideal of the diagonal $\Delta \mathcal Z$ and $I_n(\alpha_f)$ is the ideal generated by the $n\times n$ minors of $\alpha_f$ (see \cite[Proposition 3.1]{MondSomeRemarks} for the independence of $\alpha_f$). As it is the case with $K_2(f)$, the space $D^2(f)$ has a nice analytic structure whenever it has the right dimension.
\begin{prop}\label{propD2CohenMacaulay}
The dimension of $D^2(f)$ is at least $2n-p$. If $D^2(f)$ has dimension $2n-p$, then it is a Cohen-Macaulay space.
\end{prop}
As we mentioned, set theoretically, $D^2(f)$ is the projection of $K_2(f)$, that is,
\[D^2(f)=\{(u,u')\in \mathcal Y\times \mathcal Y\mid u'\neq u,f(u')=f(u)\}\cup\{(u,u)\in \Delta X\mid \operatorname{rk}(\operatorname{d}\!f_u)<n\}.\]

Going back to the reflection mapping setting, by letting $D_2^\sigma(f)$ be the image of $K^\sigma_2(f)$ in $\mathcal Y\times \mathcal Y$, we obtain a set-theoretical decomposition
\[D^2(f)=\bigcup_{\sigma\in W\setminus\{1\}}D_2^\sigma(f),\]
and clearly $D_2^\sigma(f)=\{(u,\sigma u)\in\mathcal V\times\mathcal V\mid  L(u)=L(\sigma u)=0, \ker \alpha_L^\sigma(u)\neq 0\}.$ We may improve this description, by giving  $D_2^\sigma(f)$ an adequate analytic structure:
\begin{defn}\label{defDSigma}Let $\sigma\in W\setminus\{1\}$ and $r=\dim (\Fix \sigma)$.
 If $r\geq n$, then let
		\[D_2^\sigma(f)=\{(u,\sigma u)\in\mathcal V\times\mathcal V\mid  L(u)=L(\sigma u)=0, I_{p-r}(\alpha_L^\sigma(u))=0\},\]
where $I_{p-r}(\alpha_L^\sigma)$ are the $(p-r)\times (p-r)$ minors of $\alpha_L^\sigma$. If $r< n$, then let
		\[D_2^\sigma(f)=\{(u,\sigma u)\in\mathcal V\times\mathcal V\mid  L(u)=L(\sigma u)=0\}.\]
\end{defn}
The proof that this definition is independent of the matrix $\alpha_L^\sigma$ is, mutatis mutandis, the one of \cite[Proposition 3.1]{MondSomeRemarks}.

\begin{prop}\label{propD2AnalyticStructure} For every $\sigma\in W\setminus\{1\}$, the complex space $D_2^\sigma(f)$ defined above is, as a set, the image of $K_2^\sigma(f)$ on $\mathcal Y\times \mathcal Y$, and it satisfies the following properties:
\begin{enumerate}
\item \label{propD2AnalyticStructureItem1} The dimension of $D_2^\sigma(f)$ is at least $2n-p$. If $\dim D_2^\sigma(f)=2n-p$, then $D_2^\sigma(f)$  is Cohen-Macaulay. If furthermore $\sigma$ is a reflection or $\dim \Fix \sigma<n$, then $D_2^\sigma(f)$ is locally a complete intersection. \item\label{propD2AnalyticStructureItem2} Away from $\Delta \mathcal Y$, we have an equality of complex spaces 
 \[D_2^\sigma(f)\setminus \Delta \mathcal Y=\{(u,\sigma u)\mid u\in (\mathcal Y\cap \sigma^{-1}\mathcal Y)\setminus\Fix\sigma\}.\]
\item \label{propD2AnalyticStructureItem3}Let $(u,\sigma u)\in D_2^\sigma(f)$,  assume that $u\notin \Fix\sigma$ and that $(u,\sigma u)\notin D_2^\tau(f)$, for all $\tau\neq \sigma$. Then, $D^2(f)$ and $D_2^\sigma(f)$ are locally isomorphic at $(u,\sigma u)$.

\item \label{propD2AnalyticStructureItem4}The involution $ \iota$ takes $D_2^\sigma(f)$ isomorphically to $D_2^{(\sigma^{-1})}(f)$.
\end{enumerate}
\begin{proof}We know that the image of $K_2^\sigma(f)$ in $\mathcal Y\times \mathcal Y$ is the set of pairs $(u,\sigma u)\in\mathcal Y\times\mathcal Y$, such that  $\ker \alpha_L^\sigma(u)\neq 0$. Letting $r=\dim \Fix \sigma$, the matrix $\alpha_L^\sigma$ has size $(p-n)\times (p-r)$. Therefore, $\alpha_L^\sigma$ must have non-trivial kernel automatically if $r<n$. Whenever $r\geq n$, the conditions  $\ker \alpha_L^\sigma(u)\neq 0$ and $I_{p-r}(\alpha_L^\sigma(u))=0$ are equivalent. This shows $D_2^\sigma(f)$ to be the image of $K_2^\sigma(f)$, as a set.

 Now we prove item (\ref{propD2AnalyticStructureItem1}). If $\dim \Fix \sigma<n$ then $D_2^\sigma(f)$ is defined by the $2(p-n)$ equations $L(u)=L(\sigma u)=0$ on the $p$-dimensional manifold $\{(u,\sigma u)\mid u\in \mathcal V\}$, hence it is a complete intersection if it has dimension $2n-p$. If $\sigma$ is a reflection, then $\ell_\sigma$ consists of a single entry and, then, we may compute $\alpha_L^\sigma$ just by dividing the entries of $L$ by $\ell_\sigma$, that is,
\[\label{dref}\alpha_L^\sigma=\left(\frac{\sigma^{-1}L_1-L_1}{\ell_\sigma},\dots, \frac{\sigma^{-1}L_{p-n}-L_{p-n}}{\ell_\sigma}\right)\]
Now $I_{p-r}(\alpha_L^\sigma)$ is generated by the entries of $\alpha_L^\sigma$. In this case, modulo the entries of $L$, the entries of $\sigma^{-1}L$ are generated by the entries of $I_{p-r}(\alpha_L^\sigma)$, and can be discarded as generators of the ideal of $D_2^\sigma$. Just as in the case of $\dim \Fix \sigma<n$, now $D_2^\sigma(f)$ is defined by  $2(p-n)$ equations and it is a complete intersection if it has dimension $2n-p$. If $\dim \Fix \sigma\geq n$, the fact that $D_2^\sigma(f)$  is a Cohen-Macaulay space whenever it has dimension $2n-p$ is a direct application of \cite[Thm 2.7, Lemma 2.3]{ConciniStricklandOntheVarietyOfComplexes}.

Item (\ref{propD2AnalyticStructureItem2}) is obvious if $r<n$. In the case of $r\geq n$ we need to show that, locally on a point $(u,\sigma u)\in D_2^\sigma(f)\setminus \Delta\mathcal Y$, the $p-r$ minors of
 $\alpha_L^\sigma$ are in the ideal generated by the entries of $L$ and $\sigma^{-1}L$. Let $A$
be the square submatrix obtained by picking the rows $j_1,\dots,j_{p-r}$ of $\alpha_L^\sigma$. Let $b$ be the vector with entries $\sigma^{-1}L_{j_1}-L_{j_1},\dots,\sigma^{-1}L_{j_{p-r}}-L_{j_{p-r}}$. Since $u\notin\Fix\sigma$, there must be some $i$ such that $\ell_{i}(u)\neq 0$. Let 
$A'$ be the matrix obtained by substitution of the $i$-th column of $A$ by $b$. By Cramer's Rule, we obtain $|A| = |A'|/\ell_i(u)$. The claim follows from the fact that $|A'|$ is clearly in the ideal generated by the entries of $L$ and $\sigma^{-1}L$.

Observe that item (\ref{propD2AnalyticStructureItem2}) can be restated as the fact that $D_2^\sigma(f)\setminus \Delta \mathcal Y$ and $K_2^\sigma(f)\setminus E$ are isomorphic. At the same time, $D^2(f)\setminus \Delta \mathcal Y$ and $K_2(f)\setminus E$ are isomorphic (this is true because both spaces are isomorphic to $(\mathcal Y\times_{\mathcal V}\mathcal Y)\setminus \Delta\mathcal Y$. This easy to see from the equations of $D^2(f)$ and from the definition of $K_2(f)$ given in \cite{NunoPenafortIterated}). From this point of view, item (\ref{propD2AnalyticStructureItem3}) is the same as item (\ref{propK2AnalyticStructureItem2}) from Proposition \ref{propK2AnalyticStructure}. Finally, item (\ref{propD2AnalyticStructureItem4}) is proven in the same way as item (\ref{propK2AnalyticStructureItem3}) from Proposition \ref{propK2AnalyticStructure}.
\end{proof}
\end{prop}

\begin{ex}Consider the map-germ
		\[f^{(d_1,\dots,d_5)}\colon(x,y,z)\mapsto (x^{d_1}, y^{d_2},z^{d_3},(x+y+z)^{d_4},(x-y+2z)^{d_5}),\]
 from Example \ref{k2c3c5}. With the notations that we used there,  we obtain 
		$D_2^{\sigma_1}(f)=D_2^{\sigma_2}(f)=\emptyset$ because, in both cases, $\alpha_{L}^{\sigma_i}$ has maximal rank. On the other hand, for $i=3,4,5$, $\dim Fix \sigma_i<3$, therefore
		\[D_2^{\sigma_i}(f)=\{(y,\sigma_i y)\mid L(y)=L(\sigma _iy)=0\}.\]
This spaces must be lines, because their equations are linear and we know a priori that $f^{(d_1,\dots,d_5)}$ is $\mathcal A$-finite \cite[Lemma 9.10 and Proposition 9.8]{Penafort-Sanchis:2016a}, which forces the double point space to be a reduced curve.
	\end{ex}

 \begin{rem}\label{remN+1}
In the case of $p=n+1$, every branch $D_2^\sigma(f)$ having dimension $n-1$ is locally a complete intersection. This is remarkable, since it is well known that there are germs $(\mathbb C^n,0)\to(\mathbb \C^{n+1},0)$  whose double point space $D^2(f)$ has dimension $n-1$ but is not complete intersection. Proposition \ref{propD2AnalyticStructure} shows that, for reflection mappings $\mathcal Y^n\to\mathcal V^{n+1}$, we can split $D^2(f)$ into the branches $D_2^\sigma(f)$, and these branches are complete intersections.
\end{rem}

\subsection*{Decomposition of the double point space $D(f)$}
If $f\colon \mathcal Y\to \mathcal Z$ is a finite holomorphic mapping between complex manifolds (not necessarily a reflection mapping), then the mapping $\pi\colon D^2(f)\to \mathcal Y$, given by $(y,y')\mapsto y$, is finite. The \emph{source double point space}  $D(f)$ is defined as the image of $\pi$, that is,
 \[D(f)=V(\mathcal F_0(\pi_*\mathcal O_{D^2(f)})).\]
 As a set, we have that
 \[D(f)=\{u\in \mathcal Y\mid \vert f^{-1}(f(u))\vert>1\}\cup\{u\in \mathcal Y\mid \operatorname{rk}(\operatorname{d}\!f_u)<\dim \mathcal Y\}.\]
 
 Now come back to the reflection mapping setting. Since we have the set theoretical decomposition $D^2(f)=\cup_{\sigma}D_2^\sigma(f)$, then the set $D(f)=\pi(D^2(f))$ must be the union of the sets $\pi(D_2^\sigma(f))$. 
 	\begin{defn}\label{defDf}
For any $\sigma\in W\setminus\{1\}$, we define $D_\sigma(f)=V(\mathcal F_0(\pi_*\mathcal O_{D_2^\sigma(f)}))$.	
	\end{defn}
We obtain the set theoretical decomposition
\[D(f)=\bigcup_{\sigma\in W\setminus\{1\}}D_\sigma(f).\]

Since the equations $u'=\sigma u$ are among the equations defining $D_2^\sigma(f)$, the spaces $\pi(D_2^\sigma(f))$ and  $D_2^\sigma(f)$ isomorphic. Hence the spaces $D_\sigma(f)$ may be described equivalently as follows: Let $r=\dim \Fix \sigma$.  If $r\geq n$, then 
\[D_\sigma(f)=\{u\in\mathcal Y\cap\mathcal \sigma^{-1} \mathcal Y\mid  I_{p-r}(\alpha_L^\sigma(u))=0\}.\]
If $r<n$, then 
\[D_\sigma(f)=\mathcal Y\cap\mathcal \sigma^{-1} \mathcal Y.\]

As  the spaces $D_2^\sigma(f)$ and $D_\sigma(f)$ are isomorphic, Proposition \ref{propD2AnalyticStructure} can be recast  as a result about $D_\sigma(f)$.

\begin{prop}\label{propDAnalyticStructure} For every $\sigma\in W\setminus\{1\}$, the complex space $D_\sigma(f)$ satisfies the following properties:
\begin{enumerate}
\item \label{propDAnalyticStructureItem1} The dimension of $D_\sigma(f)$ is at least $2n-p$. If $\dim D_\sigma(f)=2n-p$, then $D_\sigma(f)$  is Cohen-Macaulay. If furthermore $\sigma$ is a reflection or $\dim \Fix \sigma<n$, then $D_\sigma(f)$ is locally a complete intersection.
\item\label{propDAnalyticStructureItem2} As complex spaces, $D_\sigma(f)\setminus \Fix\sigma= (\mathcal Y\cap \sigma^{-1}\mathcal Y)\setminus\Fix\sigma$.
\item\label{propDAnalyticStructureItem3} Let $u\in D_\sigma(f)\setminus \Fix\sigma$ such that $u\notin D_2^\tau(f)$, for all $\tau\neq \sigma$. Then, $D(f)$ and $D_\sigma(f)$ are locally isomorphic at $u$.
\item \label{propDAnalyticStructureItem4}$D_\sigma(f)$ is isomorphic to $D_{\sigma^{-1}}(f)$ via $\sigma$.
\end{enumerate}
\end{prop}

If there is no risk of confusion, it is common to write $D$ instead of $D(f)$. Similarly, we may sometimes write $D_\sigma$ for  $D_\sigma(f)$. This notation appears, for example, in Example \ref{exD8Branches}.

\section{A formula for $D(f)$ in the hypersurface case}\label{secAFormula}
Let $f\colon \mathcal Y^n\to \mathcal Z^{n+1}$ be a finite mapping between complex manifolds (not necessarily a reflection mapping). The 0th Fitting ideal  $F_0(\pi_*\mathcal O_{D^2(f)})$ of the projection $\pi\colon D^2(f)\to \mathcal Y$ is then a principal ideal (if $\dim D^2(f)=n-1$, this is true because $D^2(f)$ is Cohen-Macaulay, see the first paragraph of Section \ref{secImage}. In the case of $\dim D^2(f)=n$, one simply gets the zero ideal). We usually write $\lambda\in \mathcal O_{\mathcal Y}$ for a generator of this ideal, so that
\[D(f)=V(\lambda).\]

If $f\colon \mathcal Y^n\to \C^{n+1}$ is a reflection mapping, then we have that \[D_\sigma(f)=V(\lambda_\sigma),\]
where, according to the description given below Definition \ref{defDf}, the functions $\lambda_\sigma$ are
\[\lambda_\sigma=
\begin{cases}
\frac{\sigma^{-1} L-L}{\ell_\sigma},\quad \text{if $\sigma$ is a reflection}\\
\sigma^{-1}L- L,\quad\text{otherwise}
\end{cases}
\]
 Observe that, even though $L$ vanishes on $\mathcal Y$, we have chosen to keep the term $L$ in the $\lambda_\sigma$ expressions. This allows for the divisibility by $\ell_\sigma$ to hold on the ambient space $\mathcal V$, making computations easier.

Before giving an explicit formula for the double point space $D(f)$, we need the following technical lemma, whose proof is in Appendix \ref{secDelayedProofs}.

\begin{lem}[Unfolding With Good Double Points]\label{lemIrreducibleComponentsD2}
Every multi-germ of reflection mapping $f\colon (\mathcal Y^n,Wy)\to\mathcal \C^{n+1}$ 
admits a $W$-unfolding $F$, such that, for all $\sigma, \tau\in W\setminus\{1\}$, \[\dim(D_\sigma(F)\cap \Fix\tau)<\dim D_\sigma(F).\]
 \end{lem}

%

\begin{thm}\label{thmAFormula}
As a complex space, the double point space $D(f)$ is the zero locus of
\[\lambda=\prod_{\sigma\in W\setminus\{1\}}\lambda_\sigma.\]
\begin{proof}
The expression $\Pi_{\sigma\in W\setminus\{1\}}\lambda_\sigma$ behaves well under $W$-unfoldings obviously, and it is well known that $D(f)$ behaves well under general unfoldings (see for example  \cite{Nuno-Ballesteros2015On-multiple-poi}). Therefore, we may assume $f$ to satisfy the conditions of Lemma \ref{lemIrreducibleComponentsD2} and to be generically one-to-one, by the generically one-to-one unfolding Lemma \ref{lemGen1To1Unfolding}. The generically one-to-one condition gives $\dim D^2(f)=n-1$, which forces $D^2(f)$ to be Cohen-Macaulay, by Proposition \ref{propD2CohenMacaulay}. From the decomposition $D^2(f)=\cup_{\sigma\in W\setminus\{1\}}D^\sigma_2(f)$ and Item (\ref{propD2AnalyticStructureItem1}) of Proposition \ref{propD2AnalyticStructure}, it follows that the branches $D^\sigma_2(f)$ are Cohen-Macaulay spaces of dimension $n-1$ as well. We claim that this, added to the fact that $f$ satisfies the conditions of Lemma \ref{lemIrreducibleComponentsD2}, implies that $D^2_\sigma(f)$ and $D^2_\tau(f)$ have no irreducible components in common, for all $\sigma\neq \tau$. This is true because points in $D^2_\sigma(f)\cap D^2_\tau(f)$ are of the form $(u,\sigma u)=(u,\tau u)$, hence satisfy $u\in \Fix\sigma^{-1}\tau$, but $\dim (D_\sigma(f)\cap \Fix \sigma^{-1}\tau)<n-1$.  
Now it follows from Item (\ref{propD2AnalyticStructureItem3}) of Proposition \ref{propD2AnalyticStructure} that every branch $D_2^\sigma(f)$ has a dense open subset on which $D^2(f)$ and $D_2^\sigma(f)$ are locally isomorphic. 
Since $D^2_\sigma(f)$ and  $D^2_\tau(f)$ have no common irreducible components and the isomorphisms $D_2^\sigma(f)\to D_\sigma(f)$ give $F_0((\pi_{\vert_{D^\sigma_2(f)}})_*\mathcal O_{D^\sigma_2(f)})=\langle \lambda_\sigma\rangle$, the result follows directly from Proposition \ref{propFittingIdealsDecomposition}.
\end{proof}
\end{thm}

\begin{ex}\label{exD8Branches}
Consider the  $D_8$-reflected graph 
\[f^{ D_8 }_1\colon (x,y)\mapsto (x^2+y^2,x^2y^2,2x+y),\]
from Examples \ref{exReflectedGraphs} and \ref{d8image}. With the notation of the Example \ref{exD2n}, for the reflections $\sigma_i, i=1,2,3,4$, the generators $\lambda_{\sigma_i}$ of the ideals of $D_{\sigma_i}$ are units, hence $D_{\sigma_i}=\emptyset$. The spaces $D_{\rho_1},D_{\rho_2}$ and $D_{\rho_3}$ are lines, given by the vanishing of 
 \[\lambda_{\rho_1}=x+3y, \qquad \lambda_{\rho_2}=4x+2y, \qquad \lambda_{\rho_3}=3x-y.\]
 This double point curves are depicted in Figure \ref{dreflections1}.
 Since $\rho_1^{-1}=\rho_3$, the branches $D_{\rho_1}$ and $D_{\rho_3}$ are glued together on a single branch $f(D_{\rho_1})=f(D_{\rho_3})$ of double points on the image of $f$. In contrast, since $\rho_2^{-1}=\rho_2$, the branch $D_{\rho_2}$ is glued to itself to produce the branch $f(D_{\rho_2})$ on the image. This forces $D_{\rho_2}$ to be $\rho_2$-symmetric (or, in other words, it forces $D_{\rho_2}$ to have $\rho_2$ as an automorphism). By contrast $D_{\rho_1}$ is not $\rho_1$-symmetric, but $\rho_1$ takes $D_{\rho_1}$ isomorphically onto $D_{\rho_3}$.

A more complex example is the mapping
\[f^{ D_8 }_2\colon (x,y)\mapsto (x^2+y^2,x^2y^2,2x^2+3xy-y^2+2x^3+8x^2y-2xy^2-2y^3),\]
also from Examples \ref{exReflectedGraphs} and \ref{d8image}. This time the branches associated to reflections are nonempty, and they are given by the vanishing of
\[\lambda_{\sigma_1}=-2(-3x - 8x^2 + 2y^2), \qquad \lambda_{\sigma_2}=-3x - 4x^2 - 3y - 14xy - 4y^2, \]
 \[\lambda_{\sigma_3}= 2(2x^2 + 3y - 2y^2), \qquad \lambda_{\sigma_4}=3(x - y + 2xy).\]
  The curves $D_{\sigma_i}$ are all regular, and each of them gets glued to itself to form a curve $f(D_{\sigma_i})$ on the image of $f$. Consequently, each $\sigma_i$ acts as an automorphism on $D_{\sigma_i}(f)$. The double point branches are depicted in Figure \ref{dreflections}. The elements $\rho_i$ have associated functions
  \[\lambda_{\rho_1}=3x^2 + 4x^3 + 6xy + 6x^2y - 3y^2 - 10xy^2, \qquad \lambda_{\rho_2}=4(x^3 + 4x^2y - xy^2 - y^3),\] \[\lambda_{\rho_3}=3x^2 + 6xy + 10x^2y - 3y^2 + 6xy^2 - 4y^3.\] 
The space $D_{\rho_1}$ consists of two branches $D_{\rho_1}^{(1)}$ and $D_{\rho_1}^{(2)}$. Since $\rho^{-1}_1=\rho_3$, the space $D_{\rho_3}$ is isomorphic to $D_{\rho_1}$, and its two branches $D_{\rho_3}^{(1)}$ and $D_{\rho_3}^{(2)}$ (depicted as dashed lines) are identified with the branches of $D_{\rho_1}$ to form the two branches $f(D_{\rho_1}^{(i)})=f(D_{\rho_3}^{(i)})$, with $i=1,2$. Finally, the space $D_{\rho_2}$ consists of three branches $D_{\rho_2}^{(i)}$, each being glued to itself to form the curve  $f(D_{\rho_2}^{(i)})$ on the image of $f$.

	\begin{center}
		\begin{figure}
		\includegraphics[scale=1.1]{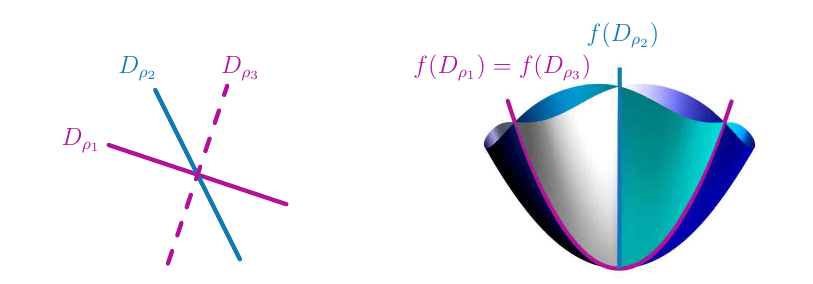}
		\caption{Double point curves of $f_1^{D_8}$}
\label{dreflections1}
	\end{figure}
	\end{center}
	\begin{center}
		\begin{figure}[h]
			\includegraphics[scale=1.1]{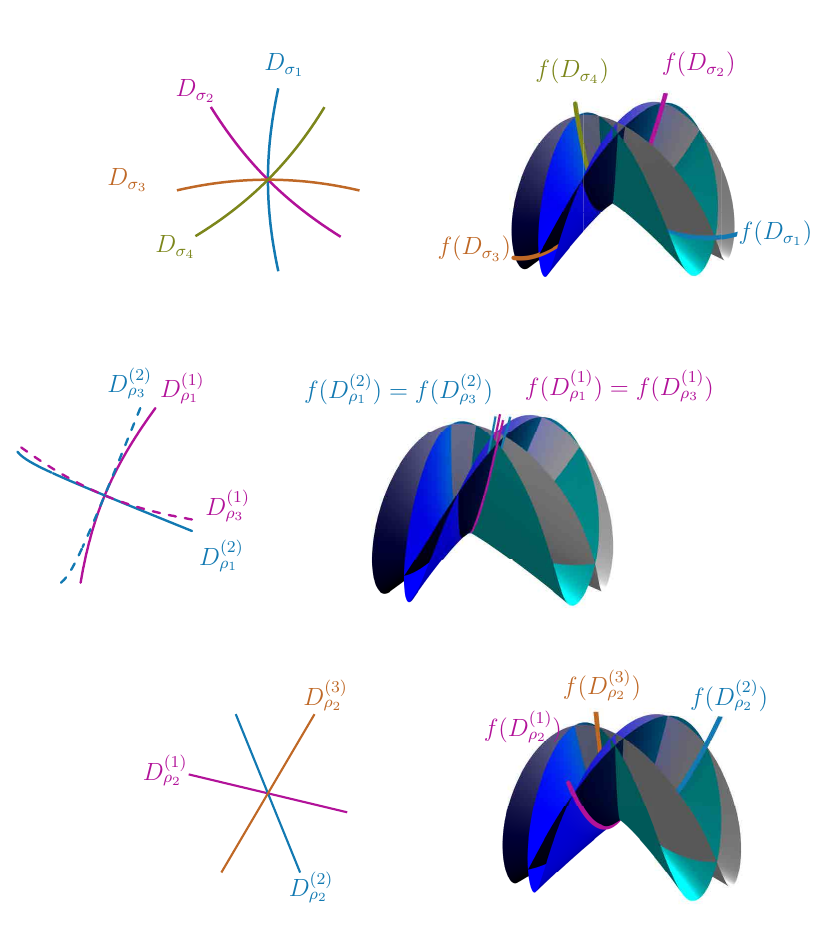}
\caption{Double point curves of $f_2^{D_8}$}
\label{dreflections}
	\end{figure}

	\end{center}


\end{ex}

\begin{rem}
It is well known that the jacobian $J\omega$ (that is, the determinant of the differential matrix of $\omega$) equals the product of the equations of the reflecting hyperplanes of all reflections in $W$. Consequently, the expression for $\lambda$ can be rewritten as
\[\lambda=\frac{\prod_{\sigma\in W\setminus\{1\}}\sigma L}{J\omega}.\]
\end{rem}

\begin{rem}\label{remDegreeOfD}
If $f=(\omega,H)$ is a reflected graph $\C^n\to\C^{n+1}$ and $H$ is homogeneous of degree $d$, as in Remark \ref{310}, then $D(f)$ is homogeneous of degree $(|W|-1)d-r$, where $r$ is the number of reflections in $|W|$.
\end{rem}

\subsection*{The double point curve $D(f)$ of a reflection mapping $\mathcal Y^2\to \C^3$}
In these dimensions the double point spaces are plane curves when they have the right dimension. Moreover, a germ $f\colon (\C^2,0)\to(\C^3,0)$ is $\mathcal A$-finite if and only if $D(f)$ is a reduced plane curve \cite{MR2868913}. For these mappings, one is interested in computing the Milnor number $\mu(D(f))$ and the delta invariant $\delta(D(f))$. This is much easier for reflection mappings than it is for arbitrary mappings, thanks to the following criterion, which follows easily by putting together Hironaka's $\mu=2\delta-r+1$ formula, the additivity of $\delta$ and the decomposition $\lambda=\Pi_{\sigma\neq 1}\lambda_\sigma$ of Theorem \ref{thmAFormula}.
\begin{prop}\label{mudelta}
A reflection mapping $f\colon (\mathcal Y^2,y)\to(\mathbb C^3,0)$ is $\mathcal A$-finite if and only if all the Milnor numbers
$\mu(D_\sigma(f))$ and all the intersection numbers $D_\sigma(f)\cdot D_\tau(f)$, with $\sigma\neq \tau$, are finite. In this case,
\[\delta(D(f))=\sum_{\sigma\in W\setminus\{1\}}\Big(\delta(D_\sigma(f))+\frac{1}{2}\sum_{\tau\in W\setminus\{1,\sigma\}}D_\sigma(f)\cdot D_\tau(f)\Big),\]

and
\[\mu(D(f))=\sum_{\sigma\in W\setminus\{1\}}\Big(\mu(D_\sigma(f))+\sum_{\tau\in W\setminus\{1,\sigma\}}D_\sigma(f)\cdot D_\tau(f)\Big)+k-|W|+2,\]
where $D_{\sigma_i}(f)\cdot D_{\sigma_j}(f)=\dim_\C(\mathcal O_\mathcal{Y}/\left\langle \lambda_{\sigma_i},\lambda_{\sigma_j}\right\rangle)$ is the intersection number of the branches $D_\sigma(f)$ and $D_\tau(f)$ and
$k$ is the number of elements $\sigma$ for which $D_\sigma(f)=\emptyset$ (We follow the convention that $\delta(\emptyset)=\mu(\emptyset)=0$).
\end{prop}

The data of the $\mu(D_\sigma(f))$, $\delta(D_\sigma(f))$, $D_{\sigma_i}(f)\cdot D_{\sigma_j}$ and whether a branch $D_\sigma(f)$ is empty can be stored easily in matrix form. Here is a convenient way to do it: First choose an ordering of the reflection group $W=\{\sigma_0,\sigma_1,\dots,\sigma_{|W|-1}\},$
where $\sigma_0=1$. Let $I$ be the $(|W|-1)\times(|W|-1)$ matrix with entries
$$
I_{ij}=\left\{\begin{matrix}
0,& \mbox{ if } i=j \mbox{ and } D_{\sigma_i}(f)\neq\emptyset,\\
1, &\mbox{ if } i=j \mbox{ and } D_{\sigma_i}(f)=\emptyset,\\
D_{\sigma_i}(f)\cdot D_{\sigma_j}(f),&\mbox{ if }i\neq j\text{ and }D_{\sigma_i}(f)\cdot D_{\sigma_j}(f)<\infty,\\
-1,&\mbox{ if } i\neq j\text{ and }D_{\sigma_i}(f)\cdot D_{\sigma_j}(f)=\infty,

\end{matrix}\right.
$$
Let $M$ be the size $(|W|-1)$ vector with entries
\[
M_i=\left\{\begin{matrix}
	0,& \mbox{ if } D_{\sigma_i}(f)=\emptyset,\\
	\mu(D_{\sigma_i}(f))&\mbox{ if } D_{\sigma_i}(f) \mbox{ is a reduced curve}\\

	-1, &\mbox{otherwise}\\
\end{matrix}\right.
\]
and define the vector $\Delta$ analogously, but replacing every $\mu(D_{\sigma_i}(f))$ with $\delta(D_{\sigma_i}(f))$. A $-1$ value anywhere indicates that $D(f)$ is non-reduced, hence that $f$ fails to be $\mathcal A$-finite.  In the abscence of any $-1$, the mapping $f$ is $\mathcal A$-finite and the formulae in \ref{mudelta} turn into
\[
\mu(D(f))=\sum_{i=1}^{|W|-1}M_{i}+\sum_{ i, j=1}^{|W|-1}I_{ij}-|W|+2,
\]
\[\delta(D(f))=\sum_{i=1}^{|W|-1}\Delta_i+\sum_{1\leq i<j\leq |W|-1}I_{ij},\]
Furthermore, the number of branches of $D(f)$ is the sum of the number of branches of all the $D_\sigma(f)$, which is equal to
\[
\sum_{i=1}^{|W|-1}(2\Delta_i-M_i-I_{ii})+|W|-1.
\]
	
Observe that the matrix $I$ is symmetric and that the vectors $M$ and $\Delta$ must contain the same value at the positions corresponding to $D_\sigma(f)$ and $D_{\sigma^{-1}}(f)$, since these are isomorphic spaces by virtue of Item \ref{propDAnalyticStructureItem4} of Theorem \ref{propDAnalyticStructureItem4}.

\begin{ex}
Consider the mapping $f_2^{D_8}$ of Example \ref{exReflectedGraphs}, whose double point branches were computed in Example \ref{exD8Branches}. By ordering the group as $1,\sigma_1,\sigma_2,\sigma_3,\sigma_4,\rho_1,\rho_2,\rho_3$ (see Example \ref{exD2n}), computing the Milnor numbers and the intersection number of all pairs of branches, one gets
\[M=(0,0,0,0,1,4,1),\quad \Delta=(0,0,0,0,1,3,1),\quad
I=\begin{pmatrix}
    0 & 1 & 1 & 1 & 2& 3 & 2\\
    1 & 0 & 1 & 1 & 2& 3 & 2\\
    1 & 1 & 0 & 1 & 2 & 3 & 2\\
    1 & 1 & 1 & 0 & 2 & 3 & 2\\
    2 & 2 & 2 & 2 & 0 & 6 & 6 \\
    3 & 3 & 3 & 3 & 6 & 0 & 6 \\
    2 & 2 & 2 & 2 & 6 & 6 & 0 \\
 \end{pmatrix}
 \] 
 From this, one computes $\mu(D(f_ 2^{D_8}))=6+104-8+2=104$ and $\delta(D(f_ 2^{D_8}))=5+52=57$.
 
This same method is efficient when applied to bigger groups. For example,  for the reflection mapping  given by the reflection group $\mathfrak S_4$ of Example \ref{examS4Group} and the embedding $(x,y)\to((x+y)^2,x,y)$, the vectors and square matrix involved have size 23, but it is better to look at the information in this format than trying to understand the space $D(f)$ formed by all branches together. For this particular example, one gets $\mu(D(f))=399$ and $\delta(D(f))=208$.
 \end{ex}

%
%

\begin{ex}
If $f$ is a reflection graph with homogeneous $H$, as in Remarks \ref{310} and \ref{remDegreeOfD}, then, by ordering $W$ as identity first, then reflections and then non-reflections, one obtains
\[M_i=\left\{\begin{matrix}	
(d-2)^2, \,  i=1,\dots,r,\\
(d-1)^2, \, i=r+1,\dots|W|-1

\end{matrix}\right.,\quad D_{\sigma}(f)\cdot D_{\tau}(f)=\left\{\begin{matrix}
(d-1)^2,  \text{if } \sigma \text{ and } \tau \text{ are reflections},\\
d^2,  \text{if } \sigma \text{ and } \tau \text{ are non-reflections},\\
d(d-1),  \text{otherwise}
\end{matrix}\right.\]
where $r$ is the number of reflections in $W$.
Therefore, \[\mu(D(f))=(1+d+r-d|W|)^2\] (observe that this also follows  from Remark \ref{remDegreeOfD}) and \[\delta(D(f))=\frac{1}{2} (-1 - d - r + d |W|) (-d - r + d |W|).\]
For instance, the map $f=f^{(d_1,d_2,d_3)}$ from Example \ref{exRefMapsUnbounded} has $\mu(D(f))=(1 - d_1 - d_2 - d_3 + d_1 d_2 d_3)^2$, and $\delta(D(f))=1/2 (1 - d_2 - d_3 + d_1 ( d_2 d_3-1)) (2 - d_2 - d_3 + d_1 (d_2 d_3-1))$.

\end{ex}

\appendix

\section{Proofs of the unfolding lemmata \ref{lemGen1To1Unfolding} and \ref{lemIrreducibleComponentsD2}}\label{secDelayedProofs}

In this section we prove the Generically One-To-One Unfolding Lemma \ref{lemGen1To1Unfolding} and the Unfolding With Good Double Points Lemma \ref{lemIrreducibleComponentsD2}, which are key in the proofs of the explicit equations of the image and double point spaces of reflection mappings in the hypersurface case.  The statements of these lemmata are quite intuitive and, if there is anything surprising about them, it is that we could not prove them more easily. One has to keep in mind however that $\mathcal Y$ is any smooth analytic of codimension one germ, with no hypothesis on how bad its relation to $W$ is, and most usual transversality arguments do not apply in this setting.  

For the most part, it is enough to apply an origin preserving rigid motion to $\mathcal Y$ in order to attain the conditions that we desire. However, certain elements of the group may impose problems that cannot be fixed in this way. We start by formalizing who the problematic elements are.

\begin{lem}Given a linear endomorphism $\sigma\colon \mathcal V\to \mathcal V$, the following statements are equivalent:
\begin{enumerate}
\item $\sigma$ preserves all $n$-dimensional vector subspaces of $V$, for certain $1<n<\dim\mathcal V$.
\item $\sigma$ preserves all vector subspaces of $\mathcal V$.
\item The matrix representation of $\sigma$ in some (or any) basis of $\mathcal V$ is
\[\sigma=
  \begin{bmatrix}
    \xi & & \\
    & \ddots & \\
    & & \xi
  \end{bmatrix}
,\]
for some $\xi\in \mathbb C$.
\end{enumerate}
\begin{proof}That the first statement implies the second follows easily after observing that, fixed $1\leq n< \dim\mathcal V$, any vector subspace can be expressed by means of intersection and sums of $n$-dimensional subspaces. The remaining implications are immediate.
\end{proof}
\end{lem}

\begin{defi}
Any linear endomorphism $\sigma\colon \mathcal V\to \mathcal V$ satisfying the conditions above is called a \emph{complex homothety} (centered at the origin). Given a subspace $K\leq \mathcal V$, we say that $\sigma$ \emph{acts as a complex homothety on $K$} if $\sigma(K)=K$ and $\sigma_{\vert_K}\colon K\to K$ is a complex homothety.
\end{defi}

We will use the following basic result, whose proof is ommited.

\begin{lem}\label{lemGrassmanian}Let $\sigma\colon \mathcal V\to \mathcal V$ be a linear automorphism and let $0\neq K\leq\mathcal V$ be a vector subspace. Let $\ell$ and $n$ be positive integers, with $\dim\mathcal V-\dim K<n<\dim \mathcal V$. Then, there is a  non-empty Zariski open subset of the product of Grassmanians,
\[\mathscr H_{\sigma,K,\ell,n}\subseteq \operatorname{Grass}(n,\mathcal V)\times \stackrel{\ell}{\dots}\times \operatorname{Grass}(n,\mathcal V),\]
such that all $H=(H_1,\dots H_\ell)\in \mathscr H_{\sigma,K,\ell,n}$ satisfy the following conditions:
\begin{enumerate}
\item $H_i\cap K\neq\sigma H_j\cap K$, for all $i$ and all $j\neq i$.
\item $H_i\cap K\neq\sigma H_i\cap K$, for all $i$, if $\sigma$ does not act as a complex homothety on $K$.
\end{enumerate}
\end{lem}

\subsection*{Proof of the Generically One-To-One Unfolding Lemma \ref{lemGen1To1Unfolding}}\label{ProofOfFirstGenericity}

Let $f$ be a multi-germ of reflection mapping, given by a germ $(\mathcal Y,Wu)$. Take the decomposition into mono-germs
\[(\mathcal Y,Wu)=(\mathcal Y_1,u^1)\sqcup\dots\sqcup (\mathcal Y_\ell,u^\ell),\]
 We want to find a representative $\mathcal Y=\mathcal Y_1\sqcup\dots \sqcup\mathcal Y_\ell$ having a trivial deformation   \[\widetilde{\mathcal Y}=\widetilde{\mathcal Y}_1\sqcup\dots\sqcup \widetilde{\mathcal Y}_\ell\subseteq \mathcal V\times \Delta,\] such that $\dim(\widetilde{\mathcal Y}\cap\sigma\widetilde{\mathcal Y})< \dim\widetilde{\mathcal Y}$, for all $\sigma\in W\setminus\{1\}$. Observe that, if the representatives $\widetilde{\mathcal Y}_i$ are chosen adequately, the desired dimension drop is equivalent to the condition that, for all $i,j$ and all $t\in \Delta\setminus\{0\}$, we have that
 \[\widetilde{\mathcal Y}_{i,t}\neq \sigma \widetilde{\mathcal Y}_{j,t}.\]
  To justify the existence of  $\widetilde{\mathcal Y}$, we consider the consider the cases $u=0$ and $u\neq 0$ separately.
  
  In the case of $u\neq 0$, we have in our advantage the fact that no complex homothety fixes any of the points in $Wu$. In particular, if $\sigma\in W$ is a complex homothety and the representatives $\mathcal Y_i$ are chosen small enough, then $\mathcal Y_i\cap \sigma \mathcal Y_i=\emptyset$. Then, any small enough deformation satisfies $\widetilde{\mathcal Y}_{i,t}\cap \sigma \widetilde{\mathcal Y}_{i,t}=\emptyset$, for all $i\leq \ell$.  
  
  By Lemma \ref{lemGrassmanian} and the Curve Selection Lemma, there exists a curve \[\gamma\colon \Delta\to \operatorname{Grass}(n,\mathcal V)\times \stackrel{\ell}{\dots}\times \operatorname{Grass}(n,\mathcal V),\] such that,  identifying $T_{u_i}\mathcal V\cong \mathcal V$ as usual, we have that $\gamma(0)=(T_{u_1}\mathcal Y,\dots, T_{u_\ell}\mathcal Y)$ and
\[\gamma(t)\in \bigcap_{\sigma\in W\setminus\{1\}}\mathscr H_{\sigma,\mathcal V,\ell,n}\] for all $t \in \Delta\setminus\{0\}$. Letting $H_{i,t}=\gamma_i(t)$, there exists an analytic family of linear automorfisms $\alpha_{i,t}\colon \mathcal V\to \mathcal V$, such that $\alpha_{i,t}(T_{u ^i}\mathcal Y)=H_{i,t}$. Finally, the germs of biholomorphism  \[\varphi_{i,t}\colon (\mathcal V,u^i)\to (\mathcal V,u^i),\] given by $u\mapsto u+\alpha_{i,t}(u-u^i)$, define deformations $\widetilde{\mathcal Y}_{i,t}=\varphi_{i,t}(\mathcal Y_i)$  of $\mathcal Y_i$, such that $u^i\in \widetilde{\mathcal Y}_{i,t}$, for all $t\in \Delta$, and $T_{u^i}\widetilde{\mathcal Y}_{i,t}=H_{i,t}$. We must check that, given $\sigma\in W\setminus\{1\}$ and assuming that either $i\neq j$ or $\sigma$ is not an homothety, the deformations $\widetilde{\mathcal Y}_i$ satisfy $\widetilde{\mathcal Y}_{i,t}\neq \sigma \widetilde{\mathcal Y}_{j,t},$ for all $t\in \Delta\setminus\{0\}$. This is trivial if $\sigma u^j\neq u^i$ and, otherwise, it follows from the fact that $\gamma(t)\in \mathscr H_{\sigma,\mathcal V,\ell,n}$, because then
\[T_{u^i}\widetilde{\mathcal {Y}}_i=H_{i,t}\neq \sigma H_{j,t}=\sigma T_{u^j}\widetilde{\mathcal {Y}}_{j,t}=T_{\sigma u^j}\sigma\widetilde{\mathcal Y}_{j,t}=T_{u^i}\sigma \tilde{\mathcal Y}_{j,t}.\]

Now consider the case where $u=0$. Here there is the advantage that the orbit of $u$ is just $\{0\}$. By the same reasoning as in the case of $u\neq 0$, there exists a family of linear automorphisms $\varphi_{t}\colon \mathcal V\to \mathcal V$, such that, $\varphi_0(\mathcal Y)=\mathcal Y$ and, for all $t\neq 0$ and  every $\sigma\in W$ which is not a complex homothety, 
\[T_0(\varphi_t(\mathcal Y))\neq T_0(\varphi_t(\sigma\mathcal Y)).\]
Now let $\sigma_1,\dots,\sigma_r$ be the homotheties in $W\setminus\{1\}$, such that $\sigma_i\mathcal Y=\mathcal Y$ (we do not need to care about  complex homotheties $\sigma$ with $\sigma \mathcal Y\neq \mathcal Y$, since they will satisfy $\widetilde{\mathcal Y}_{t}\neq \sigma \widetilde{\mathcal Y}_{t}$, for any small enough deformation of $\mathcal Y$). We may assume $\dim \mathcal Y>0$, as our problem is trivial otherwise. Then, $\mathcal Y$ cannot be contained in any $\Fix \sigma_i$, because a complex homothety which is not the identity fixes the origin only. Consequently, we may find a point $u'\in \mathcal Y$, such that $u'\neq \sigma_i(u'),$
for $i=1,\dots, r$. Since $\sigma_i(u')\neq 0$, we may take a polynomial function $p$ on $\mathcal V$, such that
 $p(u')=0$, $p(\sigma u')\neq 0$, and  $p$ vanishes at the origin with order of vanishing at least two. Now take the equations $L=(L_1,\dots,L_{p-n})=0$ of $\mathcal Y$ and define
 \[\widetilde{\mathcal Y}_t=\varphi_t\big(V(L_1+t\,p,L_2,\dots,L_{p-n})\big).\]
 
 For $t\neq 0$, the fact that $p$ vanishes at $0$ with order at least two implies that  $0\in \widetilde{\mathcal Y}_t$ and that, for every $\sigma\in W$ which is not a complex homothety, \[T_0\widetilde{\mathcal Y}_t=T_0(\varphi_t(\mathcal Y))\neq T_0(\varphi_t(\sigma\mathcal Y))=T_0(\sigma \widetilde{\mathcal Y}_t),\]
 which clearly implies $\widetilde{\mathcal Y}_t\neq \sigma \widetilde{\mathcal Y}_t$, as desired. 
 
 Finally, for  $t\neq 0$ and any of the complex homotheties $\sigma_i\in W\setminus\{\id\}$ above, it suffices to show that $\varphi_t(u')\in \widetilde{\mathcal Y}_t\setminus \sigma \widetilde{\mathcal Y}_t$. On one hand, by construction, the conditions $\varphi_t(u')\in \widetilde{\mathcal Y}_t$ and  $u'\in V(L_1+t\,p,L_2,\dots,L_{p-n})$ are equivalent, and the second holds because $u'$ is in $\mathcal Y$ and $p$ vanishes at $u'$. On the other hand, $\varphi_t(u')\in \sigma\widetilde{\mathcal Y}_t$ is equivalent to $u'\in \varphi_t^{-1}\sigma\varphi_t\big(V(L_1+t\,p,L_2,\dots,L_{p-n})\big)$.
Now observe that complex homotheties commute with linear transformations, and thus $\varphi_t^{-1}\sigma\varphi=\sigma$. Since $\sigma_i u'\in \sigma_i \mathcal Y=\mathcal Y$, we conclude $L(\sigma_iu')=0$,  and since $p(\sigma_i u')\neq 0$, it follows that $\varphi_t(u')\notin \sigma \widetilde{\mathcal Y}_t$, as desired.

\subsection*{Proof of the Unfolding With Good Double Points Lemma \ref{lemIrreducibleComponentsD2}}

Recall that our goal is to show that any multi-germ of reflection mapping admits a $W$-unfolding $F$, such that, for all $\sigma,\tau\in W\setminus\{1\}$
\[\dim(D_\sigma(F)\cap \Fix\tau)<\dim D_\sigma(F).\]

\begin{lem}\label{lemMilnorOnFacets}
Any germ of codimension one submanifold  $(\mathcal Y,Wy)\subseteq \mathcal V$ admits a $W$-unfolding, given by $\widetilde {\mathcal Y}\subseteq \mathcal V\times \Delta$, such that, for every facet $C\in \mathscr C$ and for all $t\in \Delta\setminus\{0\}$, the intersection $\widetilde{\mathcal Y}_t\cap C$ is smooth of dimension $\dim C-1$, and such that $\widetilde{\mathcal Y}$ is transverse to the reflecting hyperplanes in $\mathcal V\times \C$.
\begin{proof}
Taking the equation $L$ of $\mathcal Y$ and letting $\widetilde{\mathcal Y}=\{L(u)=t\}\subseteq \mathcal V\times \C$, the projection $\widetilde{\mathcal Y}\cap C\to \overline \C$ on the $t$ parameter is the Milnor fibration of $V(L)$ inside $\overline C$. This implies the claim that $\widetilde{\mathcal Y}_t\cap C$ is smooth of dimension $\dim C-1$. The transversality of $\widetilde{\mathcal Y}$ to the reflecting hyperplanes, which have the form $H\times \C\subseteq \mathcal V\times \C$, is obvious from its equation.
\end{proof}
\end{lem}

Since a $W$-unfolding of a $W$-unfolding of $f$ is still a $W$-unfolding of $f$ and we do not care about the number of parameters needed for the $W$-unfolding in Lemma \ref{lemIrreducibleComponentsD2}, the original $f$ may be replaced by the $W$-unfolding given by the Generically One-To-One Unfolding Lemma. This may in turn be replaced by one satisfying the conditions of Lemma \ref{lemMilnorOnFacets} as well. In other words, in order to prove Lemma \ref{lemIrreducibleComponentsD2}, we may assume our original multi-gem $f\colon (\mathcal Y^n,Wy)\to\mathcal \C^{n+1}$ to satisfy the following conditions:
\begin{enumerate}
\item \label{lemAlmostGoodUnfoldingItem1}$\mathcal Y$ intersects transversely all reflecting hyperplanes.
\item \label{lemAlmostGoodUnfoldingItem2}$\mathcal Y$ intersects properly all facets of $\mathscr C$ of dimension $n-1$.
\item \label{lemAlmostGoodUnfoldingItem3} $\dim (\mathcal Y\cap \sigma \mathcal Y)=n-1$, for all $\sigma\in W\setminus\{1\}$.
\end{enumerate}

Now observe that, to show the existence of an unfolding satisfying $\dim(D_\sigma(F)\cap \Fix\tau)<\dim D_\sigma(F)$, we do not need to deal with all pairs $\sigma$ and $\tau$ at once. We may fix $\sigma$ and $\tau$, find an a good unfolding for them and move on to the next pair of elements. We may also assume $\mathcal Y$ to have at most two branches, one containing some point $u$ and, when $\sigma u\neq u$, another one containing $\sigma u$. If more branches were present, we would solve the problem taking consecutive unfoldings, one pair of branches at a time. Additionally, condition (\ref{lemAlmostGoodUnfoldingItem2}) allows us to assume $\tau$ to be a reflection, since otherwise the dimension of $\mathcal Y\cap \Fix\tau$ is already smaller than that of $D_\sigma(f)$. By condition (\ref{lemAlmostGoodUnfoldingItem1}), the intersection $\mathcal Y\cap \Fix\tau$ is smooth of dimension $n-1$.  We are now reduced to showing the following result:

\begin{lem}
Let $\sigma\in W\setminus\{\id\}$ and let $\tau$ be a reflection in $W$. Let $\mathcal Y$ be a germ of $n$-dimensional complex manifold satisfying the conditions (\ref{lemAlmostGoodUnfoldingItem1}), (\ref{lemAlmostGoodUnfoldingItem2}) and (\ref{lemAlmostGoodUnfoldingItem3}) above, and  of one of the following forms:
\begin{itemize}
\item $\mathcal Y$ is a bigerm $(\mathcal Y_1,u)\sqcup (\mathcal Y_2,\sigma u)$ with $\sigma u\neq u$.
\item $\mathcal Y$ is a monogerm $(\mathcal Y,u)$ with $u\in \Fix \sigma$.
\end{itemize}
 Then, there exists a deformation  of $\mathcal Y$ whose fibers, for $t\neq 0$, satisfy $\dim (D_\sigma(f_t)\cap \Fix \tau^{})<n-1$.
 
 \begin{proof}
 In the bigerm case, we may proceed as in the proof of the Generically One-To-One Unfolding Lemma \ref{lemGen1To1Unfolding}, but restricting everything to $\Fix \tau$. To be precise, we start by observing that, by condition (\ref{lemAlmostGoodUnfoldingItem1}) and the fact that $W$ acts on the set of hyperplanes, the sets $\mathcal Y_{1}\cap \Fix\tau$ and $\sigma^{-1}\mathcal Y_{2}\cap \Fix\tau$ are complex manifolds of dimension $n-1$. Now  using Lemma \ref{lemGrassmanian} with $K=\Fix \tau$ and the Curve Selection Lemma, we produce a deformation of $\mathcal Y$ whose fibers $\mathcal Y_{1,t}$ and $\mathcal Y_{2,t}$ contain $u$ and $\sigma^{-1}u$, respectively, and satisfy
\[T_{u}(\mathcal Y_{1,t}\cap \Fix\tau)\neq T_{\sigma^{-1}u}(\sigma^{-1}\mathcal Y_{2,t}\cap \Fix\tau).\]
This implies $\dim (\mathcal Y_{1,t}\cap \sigma^{-1}\mathcal Y_{2,t}\cap \Fix \tau)<n-1$, and the claim follows from the inclusion $D_\sigma(f_t)\subseteq \mathcal Y_{1,t}\cap \mathcal \sigma^{-1}\mathcal Y_{2,t}$.

Now we deal with the monogerm case. Let $Z$ be the union of all $(n-1)$-dimensional components of $D_\sigma(f)\cap \Fix\tau$. Transversality forces $\mathcal Y\cap\Fix \tau$ to be a complex manifold of dimension $n-1$. Since the irreducible components of  $Z$ have dimension $n-1$ and $Z\subseteq \mathcal \Fix \tau$, we conclude that $Z$ is either empty or equal to $\mathcal Y\cap \Fix \tau$. We may assume $Z$ to be non-empty, because otherwise a trivial deformation $\mathcal Y_t=\mathcal Y$ already satisfies our claim. We distinguish two cases:

If $Z\not\subseteq \Fix\sigma$, from the inclusion $Z\subseteq D_\sigma(f)\cap \Fix\tau\subseteq \mathcal Y\cap\sigma^{-1}\mathcal Y\cap \Fix\tau$ we obtain 
\[\mathcal Y\cap \sigma^{-1}\mathcal Y\cap \Fix \tau\not\subseteq \Fix \sigma.\]
This means that, for any representative of $\mathcal Y$ (also denoted by $\mathcal Y$), there exists $u'\in \mathcal Y\cap \Fix \tau$, with $\sigma u'\neq u'$, such that $\sigma u'\in \mathcal Y$. Now take an equation $L=0$ for $\mathcal Y$, and a polynomial $p$, such that $p(u')=0$ and $p(\sigma u')\neq 0$, and define
\[\mathcal Y_t=V(L+tp).\]
For all $t\neq 0$, just as in the \hyperref[{ProofOfFirstGenericity}]{proof of the generically one-to-one unfolding Lemma \ref{lemGen1To1Unfolding}}, the fiber $\mathcal Y_t$ satisfies $u'\in \mathcal Y_t \cap \Fix \tau$ and $u'\notin\sigma^{-1}\mathcal Y_t$, in particular, $u'\notin\sigma^{-1}\mathcal Y_t\cap \Fix \tau$. This implies $\dim (\mathcal Y_{t}\cap \sigma^{-1}\mathcal Y_{t}\cap \Fix \tau)<n-1$, hence $\dim (D_\sigma(f_t)\cap \Fix\tau)<n-1$, as desired.

Finally, consider the case where $\emptyset\neq Z\subseteq \Fix \sigma$. This hypothesis requires $\mathcal Y\cap \Fix\sigma\cap \Fix\tau$ to have dimension $n-1$. Since $\Fix\sigma\cap\Fix\tau$ is the closure of a facet of the complex of $W$ and $\mathcal Y$ intersects properly the facets of dimension $n-1$, it follows that $\Fix\sigma\cap\Fix\tau$ must be a reflecting hyperplane. Then, $\sigma$ and $\tau$ are reflections with respect to this hyperplane and, since  $D_\sigma(f)\cap \Fix \tau$ depends only on $\mathcal Y,\sigma$ and $\tau$, we may assume the reflection group to be $W=\Z/d$.

 Then, the hypothesis $Z\neq \emptyset$  implies $D_\sigma(f)\neq \emptyset$, which forces the germ $f$ to be non-immersive. Since $\ker \operatorname{d}\!\omega_u=\Fix \sigma^\bot$, this implies that $T_y\mathcal Y$ contains $\Fix\sigma^\bot$. In particular, we may choose coordinates on $\mathcal V$ and $\mathcal Y$,  so that the expression of $f$ is
\[f(x_1,\dots,x_n)=(x_1,\dots,x_n^d,H(x)),\]
for some analytic function $H$ of the form
\[H(x)=H_0(x_1,\dots,x_{n-1},x_n^d)+x_n H_1(x_1,\dots,x_{n-1},x_n^{d})+ \dots+x_n^{d-1} H_d(x_1,\dots,x_{n-1},x_n^d),\]
with $H_0(0)=0$. Since $\mathcal Y$ is the graph $H$, it is isomorphic to its projection on $\C^{n}\times\{0\}\cong \C^n$ and, under this isomorphism, $\mathcal Y\cap \Fix\tau$ becomes the subset $\{x_n=0\}$ and $D_\sigma(f)$ becomes the zero locus of 
\[\lambda_\sigma=(1-\xi^\sigma) H_1(x_1,\dots,x_{n-1},x_n^d)+\dots+ (1-\xi^{(d-1)\sigma})x_n^{d-2} H_d(x_1,\dots,x_{n-1},x_n^d).\]
Clearly, we may perturb $H$ in a way that the perturbed function $\lambda_{\sigma,t}$ becomes not divisible by $x_n$ for $t\neq 0$, which forces $\dim (D_\sigma(f_t)\cap \Fix\tau)=n-2<n-1$, as desired.
 \end{proof}
\end{lem}

\bibliography{MyBibliography}
\bibliographystyle{plain}	

\end{document}